\documentclass[sn-mathphys-num]{sn-jnl}
\usepackage{graphicx}
\usepackage{mathtools}
\usepackage[section]{placeins}
\usepackage{comment}
\usepackage{amssymb}  
\usepackage{mathrsfs}

\theoremstyle{thmstyleone}
\newtheorem{theorem}{Theorem}
\newtheorem{lemma}{Lemma}
\newtheorem{remark}{Remark}
\newtheorem{prop}{Proposition}
\newtheorem{definition}{Definition}

\newenvironment{proofb}[1][\proofname]{%
  \begin{proof}[\textbf{#1}]%
}{%
  \end{proof}%
}

\title{Retrials in Loss Systems: Mean-Field Analysis, Erlang Bounds, and Power-of-$d$ Comparison}

\author{D. Cuesta}
\author{M. Jonckheere}
\author{B. Prabhu}
\author{O. Brun}
\affil{LAAS-CNRS, University of Toulouse, CNRS, 31400 Toulouse, France.\par{\{dcuesta, matthieu.jonckheere,  balakrishna.prabhu,   brun\}@laas.fr}}

\abstract{\unboldmath\mathversion{normal} We consider a system of $N$ parallel, homogeneous servers, each with finite capacity $C$. Jobs requiring one unit of capacity arrive according to a Poisson process with total rate $N\lambda$ and have exponential service times with rate $\mu$. We study a decentralized random--retrial  policy with $d$ attempts: upon arrival,   a job is assigned to a uniformly chosen server; if the server is full, the job is reassigned to another uniformly chosen server, for at most $d$ attempts; if all $d$ attempts fail, the job is rejected.
In the mean-field limit, we prove that the stationary blocking probability is characterized by the unique solution of a fixed-point equation involving the Erlang-B loss function, and we establish global convergence of the fluid model to this equilibrium. We also compare the mean-field stationary blocking probability with the analogous quantity in the power-of-$d$ model under the same parameters $(C,\lambda,\mu)$. For the same $d$, we show, via coupling techniques, that the blocking probability under power-of-$d$ is always less than or equal to the blocking under retrial-of-$d$;  but if we consider $d$ and $\tilde d$ with $\tilde d>d$, then the retrial-of-$\tilde d$ policy achieves a strictly lower blocking than the power-of-$d$ policy when $\rho=\lambda/\mu$ is sufficiently large. Moreover we show that for some specific cases, like $C=2$, $d=2$ and $\tilde d=3$, it occurs for all $\rho>0$, not only for $\rho$ sufficiently large. We complement the analysis with simulations that estimate the load level above which retrial-of-$\tilde d$ outperforms power-of-$d$ in terms of blocking.}


\begin{document}
\date{}
\maketitle

\section{Introduction}
\label{Introduction}
Loss systems with finite capacity arise in a broad range of service and communication settings, where a central question is how to reduce blocking while keeping routing decisions decentralized. On the one hand, state-aware policies that compare occupancies can substantially reduce blocking but require collecting and processing system state; on the other hand, purely random assignment has minimal information cost but can perform poorly under congestion. A standard alternative to full state information is to rely on limited sampling, by probing a small number of candidate servers per arrival.
 Two canonical ways to exploit such sampling are: (i) \emph{best-of-$d$} selection, as in the power-of-$d$ paradigm, which chooses the most favorable among $d$ sampled servers (e.g., the one with largest available capacity / smallest occupancy); and (ii) \emph{first-feasible} selection, which probes sampled servers sequentially and accepts the first one that is not full. The present paper focuses on the latter mechanism and provides a rigorous comparison with the former in a homogeneous loss setting.

We consider $N$ homogeneous servers, each with capacity $C$. Jobs arrive according to a Poisson process with total rate $N\lambda$, require one unit of capacity, and have independent exponential service times with common rate $\mu>0$, independent of arrivals and routing decisions. (We set $\mu=1$ to simplify notation). The routing rule is a decentralized random retrial policy with $d$ attempts: upon arrival, a server is sampled uniformly at random; if it is full, another server is sampled uniformly, up to $d$ total attempts; if all $d$ attempts hit full servers, the job is rejected. (We consider sampling with and without replacement). The policy uses only binary availability information, i.e., full vs. not full, and does not require comparing occupancies across sampled servers; it can be implemented as a sequential probe allowing at most $d$ attempts per arrival.

A closely related rerouting loss model (RIST) is studied in \cite{Martirosyan_Robert:2020}, where rerouted jobs are served at a different rate than jobs accepted on the first attempt (the idea is to model that reallocation  entails a cost).  Under additional assumptions, the authors show that multiple mean-field equilibrium points may arise. A different line studies retrial/orbit systems, where blocked jobs leave to an external orbit and retry after an exponential time; see e.g. \cite{Takahara96,Tuan19}. Under a link-independence assumption and in a slow-retrial regime, reduced-load fixed point approximations for blocking probability are derived, including variants in which only a fraction of blocked customers retry. By contrast, our model uses \emph{immediate} attempts at the arrival epoch (no orbit), with service rates independent of the attempt index.

Sequential first-feasible probing also arises in Video-on-Demand (VoD) systems, where a request can be served only by disks storing a replica and each disk offers a finite number of stream slots. Repeated random trials and related replica-selection schemes are analyzed in \cite{Guo2003VoD,Guo2008VoD}. Unlike our homogeneous model with a fixed probing budget $d$, these VoD schemes probe within a movie-specific replica set and are typically exhaustive up to the number of replicas. This dichotomy between sequential alternative routing (accept the first feasible alternative) and best-of-$d$ selection also appears in alternative-routing models for loss networks; see \cite{LuczakUpfal1999}.

Our first contribution is a fixed point characterization of the stationary blocking  probability. In the mean–field regime, if $\pi_C$ denotes the stationary fraction of full servers (so the stationary blocking probability is $\pi_C^d$), then $\pi_C$ is the unique solution of the scalar fixed-point equation \[ \pi_C \;=\; Erl_C\!\Big(\lambda\sum_{j=0}^{d-1}\pi_C^{\,j}\Big),\] where $Erl_C(x)=\tfrac{x^C/C!}{\sum_{i=0}^C x^i/i!}$ is the Erlang--B loss formula. We prove uniqueness of this fixed point and global convergence of the fluid model to the corresponding equilibrium; in addition, we establish finite-$N$ mean-field convergence results, including uniform convergence on compact time intervals to the limiting ODE and convergence of stationary distributions to the Dirac mass at the unique invariant state.

Our convergence arguments are analogous to those in \cite{Srikant}, which analyzes a mean-field model under the power-of-$d$ policy. There, $d\ge2$ servers are sampled uniformly and the job is placed at the sampled server with the largest available capacity; in the homogeneous unit-demand case, the fluid model has a unique invariant state and the stationary distributions converge to the corresponding Dirac mass.

A second contribution is a direct comparison between retrial-of-$d$ (first-feasible among sequential probes) and power-of-$d$ (best-of-$d$ among sampled servers). For the same $(C,\lambda,d)$, it is intuitive that power-of-$d$ yields lower blocking, since it places each arrival at the most favorable sampled server, whereas retrial-of-$d$ accepts at the first sampled server that can accommodate it. We make this intuition rigorous via a coupling construction on a common probability space, yielding a sample-path ordering that implies that power-of-$d$ has blocking probability no larger than retrial-of-$d$. The same coupling applies when the sampling degree varies over arrivals: if at the $n$th arrival both policies sample $d_n$ servers (with the same sequence $\{d_n\}_{n\ge1}$), then the ordering continues to hold.

Beyond this baseline comparison, we identify a compensation regime. For any $C$ and any $\tilde d>d$, there exists $\lambda_0=\lambda_0(C,d,\tilde d)$ such that, for all $\lambda>\lambda_0$, the retrial-of-$\tilde d$ policy achieves strictly lower blocking than power-of-$d$. Thus, while retrial-of-$d$ is simpler but (for the same $d$) incurs higher blocking than power-of-$d$, increasing the probing budget to $\tilde d>d$ can reverse this ordering once $\lambda$ is large enough. This larger probing budget should be interpreted as a worst-case budget rather than a per-arrival cost: retrial-of-$\tilde d$ probes sequentially and stops as soon as a non-full server is found, whereas power-of-$d$ inspects all $d$ sampled servers in order to compare their occupancies. Therefore, even when $\tilde d>d$, the average number of probes under retrial-of-$\tilde d$ can be smaller than $d$; in the numerical section we exhibit regimes in which retrial-of-$\tilde d$ simultaneously achieves both a lower blocking probability and a lower average number of probes than power-of-$d$. 

In addition, for $C=2$ we prove a stronger statement: retrial-of-$3$ achieves strictly lower blocking than power-of-$2$ for all $\lambda>0$. For other parameters, this uniform dominance fails at small $\lambda$, and simulations suggest a threshold behavior in $\lambda$ (depending on $(C,d,\tilde d)$) separating the regimes where retrial-of-$\tilde d$ outperforms power-of-$d$ and vice versa.

\subsection*{Organization}
The remainder of the paper is organized as follows.
Section~\ref{section_description} introduces the model and notation.
Section~\ref{section_uniqueness} characterizes the mean-field equilibrium via an Erlang--B fixed-point equation and proves uniqueness.
Section~\ref{section_convergence} is devoted to the convergence analysis.
In Section~\ref{section_comparison} we establish the coupling to compare retrial-of-$d$ with power-of-$d$, and we also analyze the regime $\tilde d>d$ in which retrial-of-$\tilde d$ achieves a smaller stationary blocking probability for sufficiently large $\lambda$ (and, in the specific case $C=2$, $d=2$, and $\tilde d=3$, for all $\lambda>0$).
Section~\ref{section_numerics} presents numerical results and discusses the inferred threshold $\lambda$ at which retrial-of-$\tilde d$ starts to outperform power-of-$d$. Proofs of auxiliary results are collected in Appendix~\ref{appendix}.

\section{Description of the model}\label{section_description}

We consider a loss system with $N$ parallel servers. All servers have the same capacity $C\in\mathbb{N}$, interpreted as the maximum number of jobs it can process simultaneously. Jobs arrive according to a Poisson process of total rate $N\lambda$ (so $\lambda>0$ is the per--server arrival rate). Service times are independent exponential random variables with common rate $\mu>0$, and they are independent of the arrival process and of all routing decisions.  We assume, without loss of generality, that each job occupies one unit of capacity while in service, and we set $\mu=1$ to simplify notation. It is worthwhile mentioning that the statements about existence and uniqueness of the stationary blocking probability extend to multiple job types $r$ with Poisson arrival rates $\lambda_r$ and exponential service rates $\mu_r$, provided each job uses the same units of capacity. In that case, the fixed-point identity holds with $\lambda$ replaced by the offered load $\sum_{r}\lambda_r/\mu_r$.

The routing policy that we consider is a decentralized one based only on instantaneous availability tests. When a job arrives, it is routed to a uniformly chosen server.  If the selected server does not have sufficient capacity to handle the job, the job is reassigned to another randomly selected server, with a maximum of $d$ attempts. If the job fails to enter the system after $d$ attempts, it  is rejected immediately. There are no queues: a job is admitted if and only if some chosen server upon inspection has available capacity. We call this policy  retrial-of-$d$.

Servers may be sampled with or without replacement (at finite $N$ they differ by hypergeometric vs.\ binomial sampling). For simplicity, we work with replacement and later show that both variants have the same fluid limit.

For $i=0,1,\dots, C$, let $S_{i}^{N}(t)$ be the fraction of servers with at least $i$ jobs in service at time $t$, and let  $\mathbf{S}^{(N)}(t)=\left\{S_{i}^{N}(t)\right\}_{i=0}^{C}$.   Note that $S_0^N(t)\equiv 1$, and define $S_{C+1}^N(t)\equiv 0$.
Let $e_i$ denote the $i$-th canonical vector in $\mathbb{R}^{C+1}$. The process $\mathbf{S}^{(N)}(t)$ is a continuous-time Markov chain, and the jumps with positive rate are the following:

\begin{align}
    &\mathbf{S}^{(N)}(t)\rightarrow \mathbf{S}^{(N)}(t)-\frac{1}{N}e_{i},\ 1\leq i\leq C,\label{initial1}\\
    &\mathbf{S}^{(N)}(t)\rightarrow \mathbf{S}^{(N)}(t)+\frac{1}{N}e_{i},\ 1\leq i\leq C.\label{initial2}
\end{align}
Transitions of type \eqref{initial1} occur when a job exits a server with exactly $i$ jobs.
Since there are $N\big(S_{i}^{N}(t)-S_{i+1}^{N}(t)\big)$ servers with exactly $i$ jobs, each departing at rate $ i$ (recall we assume, without loss of generality, $\mu=1$), this rate is equal to
\[
iN\big(S_{i}^{N}(t)-S_{i+1}^{N}(t)\big).
\]
Transitions of type  \eqref{initial2} occur when a job arrives at a server with exactly $i-1$ jobs. In this case, a job can either enter directly or after first being redirected to a full server $j$ times, with $j\leq d-1$. With replacement, the rate is
\[
\lambda N \big(S_{i-1}^{N}(t)-S_{i}^{N}(t)\big)\sum\limits_{j=0}^{d-1} \big(S_{C}^{N}(t)\big)^j.
\]
Without replacement, the rate is
\[
\lambda N\big(S_{i-1}^{N}(t)-S_{i}^{N}(t)\big)\sum\limits_{j=0}^{d-1} \frac{\binom{N S_{C}^{N}(t)}{j}}{\binom{N}{j}}\frac{N}{N-j}.
\]
This is because the first $j$ probes are directed towards full servers with probability $\frac{\binom{N S_C^N}{j}}{\binom{N}{j}}$, while the $(j+1)$-th is non-full with probability $\frac{N}{N-j}(S_{i-1}^N-S_i^N)$.

Given the transition rates of our Markov process, we introduce, as in  \cite{Srikant}, the definition of the fluid model.

\begin{definition}
Given any initial condition $\mathbf{s}^{0} \in \mathcal{S}$ with
\[
\mathcal S=
\{s\in [0,1]^{C+1}:1=s_0\geq s_1\geq\cdots\geq s_C\geq 0\},
\]
a function $\mathbf{s}(t):[0, \infty) \rightarrow \mathcal{S}$ is said to be a solution to the fluid model with initial condition $\mathbf{s}^{0}$  if: $\mathbf{s}(0)=\mathbf{s}^{0}$  and $\mathbf{s}(t)$ satisfies the following differential equations for any $t \geq 0$:
\begin{equation}\label{eq.intro.4}
    \dot{s_k}(t)=\lambda\left(s_{k-1}(t)-s_{k}(t)\right)\Big(\sum_{j=0}^{d-1}(s_C(t))^j\Big)-k\left(s_{k}(t)-s_{k+1}(t)\right),\  1 \leq k \leq C,
\end{equation}
where $s_{C+1}(t)$ is set to 0 for all $t\geq 0$.
\end{definition}
Equivalently, \eqref{eq.intro.4} can be expressed as $\dot{\mathbf{s}}(t)=\mathbf{F}\left(\mathbf{s}({t})\right)$, where
\begin{equation*}
    F_{k}(\mathbf{s})= \lambda\left(s_{k-1}-s_{k}\right)\Big(\sum_{j=0}^{d-1}s_C^j\Big)-k\left(s_{k}-s_{k+1}\right),\  1 \leq k \leq C,
\end{equation*}
and $s_{C+1}\equiv 0$.
The stationary points of the differential equation \eqref{eq.intro.4}, denoted by $\boldsymbol{\pi}$, satisfy
\begin{equation}\label{fixedpointretrial}
k\bigl(\pi_{k}-\pi_{k+1}\bigr)=\lambda\bigl(\pi_{k-1}-\pi_{k}\bigr)\Big(\sum_{j=0}^{d-1}\pi_C^j\Big),\  1 \leq k \leq C,
    \quad \pi_0=1,
\end{equation}
where $\pi_{C+1}$ is set to 0, which can also be expressed as
\begin{equation}\label{fixedpoint}
\mathbf{F}(\boldsymbol{\pi})=\mathbf{0}.
\end{equation}
We prove uniqueness of this fixed point in the next section.

We use the following notation throughout: for $C\in\mathbb{N}$ and $x\ge0$,
\[
Erl_C(x)= \frac{x^{C}/C!}{\sum_{k=0}^{C} x^{k}/k!},
\]
and
\[
\frac{\partial Erl_C(x)}{\partial x}=Erl_C'(x).
\]
\section{Uniqueness}\label{section_uniqueness}

At any stationary point $\boldsymbol{\pi}$, equation \eqref{eq.intro.4} describes an M/M/C/C queue with offered load
\[
\lambda\sum_{j=0}^{d-1}\pi_C^j,
\]
where $\pi_C$ is the stationary probability that the server is full. Therefore we have the following lemma, whose proof is included for completeness in Appendix \ref{appendix_uniqueness}.
\begin{lemma}\label{Erlangfixedpoint}
Let $\boldsymbol{\pi}$ be a solution of \eqref{fixedpoint}, then $\pi_C$ satisfies
\begin{equation}\label{eq.intro.1}
Erl_C\big(\lambda\sum_{j=0}^{d-1}\pi_C^j\big)=\pi_C.
\end{equation}
Moreover, the vector $\boldsymbol{\pi}$ is univocally determined by $\pi_C$.
\end{lemma}

We have the following result about uniqueness of the fixed point.

\begin{prop}\label{lemmaA}
There exists a unique solution $\boldsymbol{\pi} \in \mathcal{S}$ of the differential equation \eqref{eq.intro.4} that is invariant with respect to $t$, i.e., there exists a unique solution of $\mathbf{F}(\boldsymbol{\pi})=\mathbf{0}$.
\end{prop}

To prove Proposition \ref{lemmaA}, we first state an auxiliary result on the Erlang--B loss formula; its proof is deferred to Appendix \ref{appendix_uniqueness}.

\begin{lemma}\label{lemmader}
For all $C\in  \mathbb{N}$ and $x>0$,
\[
Erl_C'(x)= Erl_C(x)\left[\frac{C}{x}-1+ Erl_C(x)\right ]<\frac{1}{x}\left[1- Erl_C(x)\right].
\]
\end{lemma}

Now we present the proof of Proposition \ref{lemmaA}.
\begin{proofb}[Proof of Proposition \ref{lemmaA}]
Let $g:[0,1]\rightarrow[0,+\infty)$ be defined by $g(p)=\sum_{j=0}^{d-1} p^j$. By Lemma \ref{Erlangfixedpoint}, $\boldsymbol{\pi}$ is univocally determined by $\pi_C$, and $\pi_C$ satisfies
\[
Erl_C(\lambda g(\pi_C))=\pi_C.
\]
So it is sufficient to prove that for all $C\in \mathbb{N}$ and $\lambda>0$, the function
\[
P_{C,\lambda}(p)= Erl_C\big(\lambda g(p)\big)-p
\]
has a unique real root in $[0,1]$.

Since $P_{C,\lambda}(0)= Erl_C(\lambda)>0$ and $P_{C,\lambda}(1)= Erl_C(\lambda d)-1<0$, we know there is at least one root in the interval $(0,1)$. To prove there is exactly one, let us first show that for any root $p_0$ of $P_{C,\lambda}(p)$, it holds that $P'_{C,\lambda}(p_0)<0$, where the derivative is with respect to $p$.

Let $p_0$ be a root of $P_{C,\lambda}(p)$. By Lemma \ref{lemmader},
\[
Erl_C'(x)<\frac{1}{x}(1- Erl_C(x)) \text{ for all } x>0,
\]
which implies
\begin{align}\label{newequation7}
P'_{C,\lambda}(p)&= Erl_C'(g(p)\lambda)g'(p)\lambda -1\nonumber\\
&<\frac{g'(p)}{g(p)}(1- Erl_C(g(p)\lambda))-1,
\end{align}
where in the last step we also used $g'(p)>0$ because $g(p)=\sum_{j=0}^{d-1} p^j$ is increasing. Moreover, $g(p)=\frac{1-p^d}{1-p}$ and
\[
g'(p)=\frac{-dp^{d-1}(1-p)+(1-p^d)}{(1-p)^2},\ p\in(0,1),
\]
and
\begin{equation}\label{newequation8}
\frac{g'(p)}{g(p)}=\frac{-dp^{d-1}}{1-p^d}+\frac{1}{1-p}<\frac{1}{1-p}.
\end{equation}
From \eqref{newequation7} and \eqref{newequation8} it follows that
\begin{align}\label{newequation9}
P'_{C,\lambda}(p)<\frac{1}{1-p}(1- Erl_C(g(p)\lambda))-1.
\end{align}
Since $p_0$ is a root of $P_{C,\lambda}(p)$, $ Erl_C(\lambda g(p_0))=p_0$, and from \eqref{newequation9} we obtain
\begin{align}\label{neweq4..}
P'_{C,\lambda}(p_0)<\frac{1}{(1-p_0)}(1-p_0)-1=0.\nonumber
\end{align}

Now let us prove uniqueness of the root by contradiction. Let us suppose there exists more than one root in the interval $(0,1)$ and let $p_1<p_2$ be the two smallest roots, with $p_1<p_2$ (there are no double roots because $P'_{C,\lambda}(p_0)<0$ for all $p_0$ root). Since $P_{C,\lambda}(0)>0$ and $P_{C,\lambda}(p)$ is continuous, then $P_{C,\lambda}(p)$ is negative in the interval $(p_1,p_2)$; in particular, it is negative in the interval $(p_2-\epsilon,p_2)$ for $\epsilon>0$ sufficiently small. But also $P_{C,\lambda}(p)$ is positive in  $(p_2,p_2+\epsilon)$ for $\epsilon>0$ sufficiently small, which contradicts $P'_{C,\lambda}(p_2)<0$.
\end{proofb}
\section{Convergence results}\label{section_convergence}

Having established uniqueness, we present a regularity statement for the fluid vector field, showing that the drift $\mathbf{F}$ in \eqref{eq.intro.4} is Lipschitz on $\mathcal S$. We then state the convergence results. Since the arguments follow \cite{Srikant} with minor modifications specific to the retrial-of-$d$ routing, the proofs are postponed to Appendix \ref{appendix_convergence}.

\begin{lemma}\label{lemmaD}
The drift function $\mathbf{F}(\cdot)$ in \eqref{eq.intro.4} is Lipschitz on $\mathcal S$, i.e., there exists a constant $M>0$ such that for any $\mathbf{x}, \mathbf{y} \in \mathcal{S}$,
\[
|\mathbf{F}(\mathbf{x})-\mathbf{F}(\mathbf{y})| \leq M|\mathbf{x}-\mathbf{y}|,
\]
where $|\cdot|$ denotes the $\ell_{1}$ norm. Moreover, viewed as a map on $\mathbb R^{C+1}$, $\mathbf F(\cdot)$ is locally Lipschitz.
\end{lemma}

\begin{lemma}\label{lemmaB}
For every initial condition $\mathbf{s}^0\in\mathcal S$, the fluid model \eqref{eq.intro.4} admits a unique solution $\mathbf{s}(\mathbf{s}^0,t)$ defined for all $t\ge0$. Moreover, the set $\mathcal S$ is forward invariant, i.e., $\mathbf{s}(\mathbf{s}^0,t)\in\mathcal S$ for all $t\ge0$.
\end{lemma}

\begin{lemma}\label{lemmaB2}
Let $\mathbf{s}(\mathbf{s}^0,t)$ be a solution for the fluid model \eqref{eq.intro.4} with initial condition $\mathbf{s}^0$ and let $\boldsymbol{\pi}$ be the unique stationary point of \eqref{eq.intro.4} given by Proposition \ref{lemmaA}. Then $\mathbf{s}(\mathbf{s}^0,t)$ converges to $\boldsymbol{\pi}$ as $t\to\infty$.
\end{lemma}

\begin{lemma}\label{lemmaC}
Consider a sequence of systems with the number of servers $N$ increasing to infinity. Fix any $T>0$. If the initial states satisfy $\mathbf{S}^{(N)}(0)\to \mathbf{s}^{0}\in\mathcal{S}$ almost surely as $N\to\infty$, then
\begin{equation*}
\lim _{N \rightarrow \infty} \sup _{t \in[0, T]}\left|\mathbf{S}^{(N)}(t)-\mathbf{s}\left(\mathbf{s}^{0}, t\right)\right|=0 \quad \text { almost surely.}
\end{equation*}
where $\mathbf{s}\left(\mathbf{s}^{0}, t\right)$ is the solution to the differential equation \eqref{eq.intro.4} with initial condition $\mathbf{s}^{0}$.
\end{lemma}

\begin{theorem}\label{theorem1}
Let $\boldsymbol{\pi}$ be the unique stationary point of \eqref{eq.intro.4} given by Proposition \ref{lemmaA}. Then the sequence of stationary distributions $\boldsymbol{\pi}^{(N)}$ of $\mathbf{S}^{(N)}$ converges weakly to $\delta_{\boldsymbol{\pi}}$, the Dirac mass at $\boldsymbol{\pi}$.
\end{theorem}

We conclude this section with two asymptotic results as $N\to\infty$: propagation of chaos, and convergence of the stationary blocking probability to $\pi_C^{\,d}$.

\begin{prop}[Propagation of chaos]\label{prop_chaos}
Let ${\Pi}^{(N)}$ be the joint stationary distribution of the $N$-server system on $E^N$, with $E=\{0,1,\dots,C\}$, and let ${\Pi}^{(N)}_{\mid l}$ denote its marginal on the first $l$ coordinates, for fixed $l\in\mathbb{N}$. Let $\boldsymbol{\pi}=(\pi_0,\dots,\pi_C)$ be the unique stationary point of \eqref{eq.intro.4}, and define $\rho$ the probability measures on $E$ given by
\begin{equation}\label{eqrho}
\rho(k)= \pi_k-\pi_{k+1},\qquad 0\le k\le C,\quad (\pi_{C+1}= 0).
\end{equation}
 Then, as $N\to\infty$, the sequence of measures ${\Pi}^{(N)}_{\mid l}$ converges weakly in $\mathcal P(E^l)$ to the product measure $\rho^{\otimes l}$, where $\rho^{\otimes l}$ is the joint law of $l$ independent copies of $\rho$.
\end{prop}

\begin{prop}\label{prop_limit_blocking}
Let $P_{block}^{(N)}$ be the stationary blocking probability experienced by arriving jobs in a system with $N$ servers. Then $P_{block}^{(N)}$ converges to $\pi_C^d$ as $N$ tends to $+\infty$.
\end{prop}

\section{Comparison against power-of-d policy} \label{section_comparison}
We recall the power-of-$d$ policy. As in the retrial-of-$d$ policy, we consider $N$ servers of capacity $C$. Jobs arrive at total rate $N\lambda$ (so $\lambda$ is the per–server rate). Each job occupies one unit of capacity. Upon arrival, $d$ servers are sampled uniformly; the job is routed to the least-loaded among them (ties are broken uniformly). Service times are exponential with rate $\mu$; and as before, we set $\mu=1$ for notational convenience. Under this policy, the mean-field dynamics has a unique fixed point $\boldsymbol{q}$ (see, e.g., \cite{Srikant}) satisfying
\begin{equation}\label{eqpower}
  k\bigl(q_{k}-q_{k+1}\bigr)=\lambda\bigl(q_{k-1}^d-q_{k}^d\bigr),
\ 1 \leq k \leq C,
    \quad q_0=1,\end{equation} 
where $q_{C+1}$ is set to 0.

Given $C\in\mathbb{N}$ and $\lambda>0$, it is intuitive that power-of-$d$ yields smaller or equal blocking than retrial-of-$d$. (Indeed, retrial-of-$d$  can be interpreted as: sample $d$ servers uniformly; if all $d$ are full, block; otherwise, route uniformly  among the sampled non-full servers).  In the next subsection we present a formal coupling to prove that. 

\subsection{Coupling}
Let us consider the two systems, power-of-$d$ and retrial-of-$d$, with the same parameters $C\in \mathbb N$, $\lambda>0$, $\mu>0$, and the same number of servers $N\in \mathbb{N}$. Both systems share the same arrival and service parameters: arrivals follow a Poisson process of rate $N\lambda$, and each job has an exponential service time with rate $\mu$.

In each system, in addition to the $N$ real servers (with indices $1$ through $N$), we introduce an extra ``fictitious'' server indexed by $0$. Rejected jobs are deposited in server $0$ in the corresponding system.
 
For each $0 \le i \le N$, let $p_i$ and $r_i$ denote the occupancy at server $i$ in the power-of-$d$ and  the retrial-of-$d$ systems, respectively. For $1 \le i \le N$, both $p_i$ and $r_i$ are upper bounded by the total capacity $C$, whereas $p_0$ and $r_0$ are unbounded. Moreover, in both systems, server $0$ has no departures. Thus, $p_0$ and $r_0$ equal the cumulative number of rejected jobs in the power-of-$d$ and the retrial-of-$d$ systems, respectively.

In the following, we denote by $\mathscr{S}_P$  the system operated under the power-of-$d$ scheme, and by $\mathscr{S}_R$ the system operated under the retrial-of-$d$ scheme. Also we set $p=(p_0,p_1,\dots,p_N)$ and $r=(r_0,r_1,\dots,r_N)$.  We now describe the coupling.

\textbf{Arrivals:} since the arrival rate is the same in both systems ($N\lambda$), we use the same Poisson arrival times, so arrivals occur simultaneously in both systems. At each arrival, we sample $d$ indices uniformly from $\{1,\dots,N\}$. We denote them $i_1,\dots,i_d,$ with
$$1\le i_1\le i_2\le \cdots \le i_d\le N.$$
(With replacement, the inequalities between consecutive indices are non-strict; without replacement, they are strict. The coupling works in both variants.)
\begin{itemize}
\item Retrial-of-$d$ ($\mathscr{S}_R$): among the $d$ sampled servers that are not full (occupancy less than $C$), select one uniformly at random. If all the sampled servers are full, select the fictitious server $0$. The selected server receives the job.
\item Power-of-$d$ ($\mathscr{S}_P$): among the $d$ sampled servers, select one with minimum occupancy (ties are broken uniformly at random). If all the sampled servers are full, select the fictitious server $0$. The selected server receives the job.
\end{itemize}
Separately in each system ($\mathscr{S}_P$  and $\mathscr{S}_R$), immediately after placing the job, we apply the following \textit{relabeling procedure} (which is also applied after each departure of a job). For $\mathscr{S}_P$, we relabel the indices of servers $1$ through $N$ so that the occupancies are in non-increasing order:
\[
p_1 \ge p_2 \ge \cdots \ge p_N.
\]
This is purely a relabeling and does not affect the dynamics, since sampling $d$ servers uniformly at random is invariant under relabeling.
We then label the jobs in servers $1$ through $N$ with labels $1,2,\dots,J_p$, where $J_p = p_1 + \cdots + p_N$, ordered first by server index and, within each server, arbitrarily. Concretely, jobs at server $1$ receive labels $1$ to $p_1$, jobs at server $2$ receive labels $p_1+1$ to $p_1+p_2$, and so on.  (Notice that the fictitious server $0$ and the jobs it contains are excluded from this relabeling procedure). We repeat the same procedure for $\mathscr{S}_R$: we relabel the indices of servers $1$ through $N$ so that
\[
r_1 \ge r_2 \ge \cdots \ge r_N,
\]
and we label the jobs from 1 to  $J_r = r_1 + \cdots + r_N$ such that jobs at server $1$ receive labels $1$ to $r_1$, jobs at server $2$ receive labels $r_1+1$ to $r_1+r_2$, and so on.  

\textbf{Departures:} Suppose that in both systems, the occupancies of servers $1$ through $N$ are non-increasing, and jobs are labeled in such a way that server $1$ carries the smallest labels, server $2$ the next block of labels, and so on. Since all jobs have the same exponential rate $\mu$, we couple the departure clocks as follows. 
Let $J_p$ and $J_r$ be the current numbers of jobs in servers $1$ to $N$ in $\mathscr{S}_P$  and $\mathscr{S}_R$ respectively, and let $J=\max\{J_p,J_r\}$. We run $J$ independent exponential clocks of rate $\mu$ numbered from $1$ to $J$, and let $n$ be the index of the clock that rings. Then the job with label $n$ (if any) leaves the system (it can happen in both systems $\mathscr{S}_P$  and $\mathscr{S}_R$  or only in one). This coupling is consistent with the departure rules in both models. As we did before, after each departure we apply the same relabeling procedure separately in $\mathscr{S}_P$  and $\mathscr{S}_R$.

Notice that, since we apply the relabeling procedure immediately after each arrival and also each departure, at all times the occupancies of servers $1$ through $N$ are non-increasing, and jobs are labeled in such a way that server $1$ carries the smallest labels, server $2$ the next block of labels, and so on.

\textbf{Equivalent interpretation of arrivals and departures:} Because servers are relabeled after each arrival and each departure, we may equivalently view placement as follows. If the job is rejected, the job is added to the fictitious server $0$. Otherwise, in $\mathscr{S}_R$, we add the job to the smallest-index server among servers $1$ through $N$ whose occupancy equals that of the chosen server; in $\mathscr{S}_P$, we add the job to the smallest-index server among servers $1$ through $N$ whose occupancy equals that of the server with index $i_d$.

Likewise, we may view departures as follows. If the job with label $n$ is in a server with $t$ jobs, we remove one job from the largest-index server among those with occupancy $t$. 

This new and equivalent way to see arrivals and departures will be more useful, so from now on, when we refer to arrivals and departures, we are directly referring to this equivalent view.

We now introduce two useful definitions for comparing the sequences associated with both systems.  
\begin{definition}
We say a sequence $a=(a_0,a_1,\dots,a_N)$ is admissible if $a_i\in \mathbb{N}_0$ for all  $i=0,1,\cdots, N$ and $C\geq a_1\geq a_2\geq  \cdots\geq  a_N$.  \end{definition}

Notice that if $a$ and $b$ are admissible sequences, then after any \textit{arrival} or \textit{departure}, the resulting sequences remain admissible.

\begin{definition}
    Given any sequence $a=(a_0,a_1,\dots,a_N)$, let $$S_i(a)=\sum_{k=0}^i a_k$$ be the sum of the first $i+1$ coordinates. We say that the sequence $a=(a_0,a_1,\dots,a_N)$ majorizes by sums the sequence $b=(b_0,b_1,\dots,b_N)$ (and we denote this by $a \succeq_{\Sigma} b$) if both sequences are admissible and
$$ S_i(a)\geq S_i(b)  \text{ for all }i=0,1,\dots,N.$$ 
\end{definition}

In the next two lemmas, we prove that if a sequence majorizes by sums another sequence, then after any \textit{arrival} or \textit{departure}, the majorization-by-sums property is preserved.
 
\begin{lemma}\label{coupling1_arr}
Fix $t\geq 0$ and let $r=(r_0,r_1,\dots,r_N)$ and $p=(p_0,p_1,\dots,p_N)$ be the sequences associated to the retrial-of-$d$ and power-of-$d$ systems at time $t$ under the coupling. Let $r'$ and $p'$ be the sequences obtained from $r$ and $p$ respectively after an \textit{arrival}. If $r \succeq_{\Sigma} p$, then $r' \succeq_{\Sigma} p'$.
\end{lemma}

\begin{proofb} 
 Suppose that $r\succeq_{\Sigma}p$ and that after an \textit{arrival}, the new sequence $r'$ does not majorize by sums the new sequence $p'$. We know that each \textit{arrival} increases exactly one job in each system (recall that, if all the $d$ selected servers in a system are full, then the fictitious server with index 0 increases by 1). Let $i$ be the index of the server in $\mathscr{S}_P$ to which the job is added, and $l$ be the index  in $\mathscr{S}_R$ to which the job is added. Then, if $i\geq l$, or $i<l$ and $S_k(r)>S_k(p)$ for all $k\in\{i,\dots,l-1\}$, it holds that $r'\succeq_{\Sigma} p'$. Hence, we assume that there exists $j\in\{0,\dots,N-1\}$ such that the job in $\mathscr{S}_P$  was added to a server with index in $\{0,1,\dots,j\}$, the job in $\mathscr{S}_R$  was added to a server with index in $\{j+1,\dots,N\}$, and $S_j(r)=S_j(p)$. We consider the largest $j$ with that property, that is, the largest $j\in\{0,\dots,N-1\}$ such that $i\in\{0,1,\dots,j\}$, $l\in \{j+1,\dots,N\}$ and $S_j(r)=S_j(p)$, and show that this cannot happen.

First, we notice that from $S_{j+1}(r)\geq S_{j+1}(p)$ and $S_{j}(r)= S_{j}(p)$, it follows $r_{j+1}\geq p_{j+1}$. 
Since the job in $\mathscr{S}_R$  was added to a server with index at least $j+1$, the fictitious server with index 0 did not increase by 1. Hence the server that increased by 1 was among the selected servers $i_1\leq \dots\leq i_d$, and therefore $j+1\leq i_d$. Since in $\mathscr{S}_P$  the job was added to a server with index at most $j$, the job could not be added to the server with index $i_d$. Moreover, in $\mathscr{S}_P$, if all the selected servers were full (equivalently, if $p_{i_d}=C$), then the job was added to the fictitious server with index 0; otherwise, it was added to the smallest-index server among servers 1 through $N$ whose occupancy equaled that of server $i_d$. Consequently, either $p_{i_d}=C$ and hence, $p_1=\cdots=p_{i_d}=C$, or $p_{i_d}<C$ and hence, $j\geq 1$ and $ p_j=\cdots=p_{i_d}.$

In the first case, from  $C\geq r_{j+1} \geq p_{j+1} = C$ we obtain  $C=r_{j+1}=p_{j+1}$ and, combined with the fact that $S_j(r)=S_j(p)$, it follows $S_{j+1}(r)=S_{j+1}(p)$. Since $S_{j+1}(r)=S_{j+1}(p)$ and the added job in $\mathscr{S}_R$ was to a server with index at least $j+1$, it had to be added to the server with index $j+1$; otherwise, the maximality of $j$ would be broken ($j+1$ would satisfy the same property that $j$). But it is not possible, because $r_{j+1}=C$ and hence that server was full.  

In the second case, we have $r_j\geq r_{j+1}\geq p_{j+1}=p_j$. Since $S_{j-1}(r)\geq S_{j-1}(p)$ and $S_{j}(r)=S_{j}(p)$, it follows that $r_j\leq p_j$. Then $r_j = r_{j+1}= p_{j+1} = p_{j}$, and therefore $S_{j-1}(r)=S_{j-1}(p)$ and $S_{j+1}(r)=S_{j+1}(p)$. Moreover, since $r_j=r_{j+1}$, the job in $\mathscr{S}_R$  could not be added to server $j+1$ (because there was a server with a smaller index having the same number of jobs), and since it was added to a server with index at least $j+1$, it had to be added to a server with index at least $j+2$; but again this would break the maximality of $j$ (because $j+1$ would satisfy the same property). Therefore, we conclude that after any \textit{arrival},  $r'\succeq_{\Sigma} p'$. \end{proofb}

\begin{lemma}\label{coupling1_dep}
Fix $t\geq 0$ and let $r=(r_0,r_1,\dots,r_N)$ and $p=(p_0,p_1,\dots,p_N)$ be the sequences associated to the retrial-of-$d$ and power-of-$d$ systems at time $t$ under the coupling. Let $r'$ and $p'$ be the sequences obtained from $r$ and $p$ respectively after a \textit{departure}. If $r \succeq_{\Sigma} p$, then $r' \succeq_{\Sigma} p'$.
\end{lemma}

\begin{proofb} 
 
 Suppose that $r\succeq_{\Sigma} p$ and after a \textit{departure}, the new sequence $r'$ does not majorize by sums the new sequence $p'$. Each \textit{departure} removes at most one job in each system. Moreover, if after a \textit{departure}: 
\begin{itemize}
    \item a job was removed in $\mathscr{S}_P$  from a server $i$ and in $\mathscr{S}_R$  from a server $l$ with $i\leq l$, or a job was removed only in $\mathscr{S}_P$;
\item or a job was removed in $\mathscr{S}_P$  from a server $i$ and in $\mathscr{S}_R$  from a server $l$ with $i>l$ and $S_k(r)>S_k(p)$ for all $k\in\{l,\dots,i-1\}$;
\item or a job is only removed in $\mathscr{S}_R$  from a server $l$ and $S_k(r)>S_k(p)$ for all $k\in\{l,\dots,N\}$,
\end{itemize} then $r'\succeq_{\Sigma} p'$. Hence, we assume that there exists $j\in\{1,\dots,N\}$ such that a job in the $\mathscr{S}_R$  was removed from one of the servers $\{1,\dots,j\}$, no job on servers $\{1,\dots,j\}$ was removed in $\mathscr{S}_P$, and $S_j(r)=S_j(p)$. We consider the smallest $j$ with that property, and show that this cannot happen.

Let $n$ be the number of the departure clock that rang (recall that jobs in the fictitious server $0$ have no associated labels or clocks). Since the removed job in $\mathscr{S}_R$  was at a server with index in $\{1,\cdots,j\}$, the job with label $n$ in $\mathscr{S}_R$  was at a server with index in $\{1,\cdots,j\}$, which means $n\leq S_j(r)-r_0$. Since  $r_0=S_0(r)\geq S_0(p)= p_0$, $n+p_0\leq n+r_0 \leq S_j(r)=S_j(p)$. Then $n\leq S_j(p)-p_0$ and therefore, the job with label $n$ in $\mathscr{S}_P$  was at a server with index in $\{1,\cdots,j\}$ as well. Since there is a job with label $n$ in $\mathscr{S}_P$, a job was removed in this system as well. Since that removal did not occur at any of the servers $\{1,\cdots,j\}$ but it occurred at a server with the same number of jobs as the server holding the job with label $n$, then we have $p_{j+1}=p_j$. Since $S_{j+1}(r)\geq S_{j+1}(p)$ and $S_{j}(r)= S_{j}(p)$, then $r_j\geq r_{j+1}\geq p_{j+1}=p_j$. Since $r_j\geq p_j$, $S_{j}(r)=S_{j}(p)$ and $S_{j-1}(r)\geq S_{j-1}(p)$, it follows that $S_{j-1}(r)= S_{j-1}(p)$ and $r_j=p_j$. Since $S_{j-1}(r)= S_{j-1}(p)$, the removed job in $\mathscr{S}_R$ had to be removed from server $j$, otherwise the minimality of $j$ would be broken ($j-1$ would satisfy the same property that $j$). Since the job was removed from server $j$, it follows that $j$ was the largest index among  servers 1 through $N$ that had the same number of jobs as the server holding the job with label $n$. But it is false because $r_j=r_{j+1}$. Therefore, we conclude that after any \textit{departure},  $r'\succeq_{\Sigma} p'$. \end{proofb}

\begin{prop}
    \label{coupling2} Fix $N\in \mathbb{N}$ and let $$r^{(N)}(t)=(r^{(N)}_0(t),r^{(N)}_1(t),\dots,r^{(N)}_N(t)),\qquad
p^{(N)}(t)=(p^{(N)}_0(t),p^{(N)}_1(t),\dots,p^{(N)}_N(t))$$
be the reordered occupancy vectors in retrial-of-$d$ and power-of-$d$ respectively, at time $t\geq 0$, under the coupling. Then, for all $t\ge 0$,
$$r^{(N)}(t)\  \succeq_{\Sigma}\ p^{(N)}(t).$$
In particular, the cumulative fraction of rejected jobs in retrial-of-$d$ at time $t$ is greater than or equal to the cumulative fraction of rejected jobs in the power-of-$d$ at time $t$ for all $t\geq 0$.
\end{prop}

\begin{proofb}  
Since $r^{(N)}(t)$ and $p^{(N)}(t)$ are obtained from $r^{(N)}(0)$ and $p^{(N)}(0)$ after a countable sequence of \textit{arrival}s and \textit{departure}s, and at time $0$  $$r^{(N)}(0)\  \succeq_{\Sigma}\ p^{(N)}(0)$$ because all the servers start empty,  by Lemma \ref{coupling1_arr} and Lemma \ref{coupling1_dep} it follows $$r^{(N)}(t)\  \succeq_{\Sigma}\ p^{(N)}(t).$$
Taking $k=0$ yields $r^{(N)}_0(t)\ge p^{(N)}_0(t)$, and dividing by the common number of arrivals up to time $t$, we obtain the inequality for the rejection fractions.
\end{proofb}

\begin{theorem}
    \label{coupling3}
Let $\pi=(\pi_0,\pi_1,\dots,\pi_C)$ and $q=(q_0,q_1,\dots,q_C)$ denote the unique mean-field stationary solutions of the retrial-of-$d$ and power-of-$d$ given by \eqref{fixedpointretrial} and \eqref{eqpower} respectively. Then
$$(\pi_C)^d\ \geq\ (q_C)^d.$$
\end{theorem}

\begin{proofb} 
By Proposition \ref{coupling2}, the cumulative fraction of rejected jobs in retrial-of-$d$ at time $t$ is greater than or equal to the cumulative fraction of rejected jobs in the power-of-$d$ at time $t$ for all $t\geq 0$ and for all $N\in \mathbb{N}$.
Since these fractions converge to the stationary blocking probabilities of the $N$-systems, which we denote with $P_{\mathrm{block}}^{(N)}$ and $P_{\mathrm{block}}^{(N),\mathrm{power}}$ for retrial-of-$d$ and power-of-$d$ respectively, we obtain
$$P_{\mathrm{block}}^{(N)} \geq P_{\mathrm{block}}^{(N),\mathrm{power}}
\ \text{ for each } N\in\mathbb{N}.$$
Letting $N\to\infty$, by Proposition \ref{prop_limit_blocking} and the analogous result for the power-of-$d$ model (Section 2.1 in \cite{Srikant}), we obtain  $(\pi_C)^d \geq (q_C)^d$.  
\end{proofb}

\subsection{Comparison in high-load regime}

In this subsection, we show that in the high-load regime, adding even a single extra probe in the retrial scheme yields strictly lower blocking than the power-of-$d$ policy.
\begin{theorem}\label{theorem2}
Fix integers $1\le d<\tilde d$. There exists $\lambda_0=\lambda_0(C,d,\tilde d)$ such that for all $\lambda>\lambda_0$, the stationary blocking probability under retrial-of-$\tilde d$ is strictly smaller than under power-of-$d$.
\end{theorem}

The proof proceeds by deriving Erlang--type bounds on the stationary blocking probabilities for both policies and then comparing these bounds in the high-load regime.

Fix $C\in\mathbb{N}$. We start by deriving an Erlang--type upper bound for the stationary blocking probability $\pi_C^{d}$ under retrial-of-$d$.
\begin{prop}\label{propretrial}
Fix $d\geq 1$, $C\in \mathbb{N}$ and $\lambda> C$. Let $\boldsymbol{\pi}$ be the fixed point of equation \eqref{fixedpoint} associated with the retrial-of-$d$ scheme. Then
\begin{equation}\label{upperboundretrial}
\pi_C^{d}\leq [ Erl_C(\mathcal{R}_{{d}}(\lambda,C))]^{d}
\end{equation}
with
\[
\mathcal{R}_{d}(\lambda,C)=\frac{C(1-\tfrac{C}{\lambda})^{1/{d}}}{1-\Bigl[1-\frac{C}{\lambda}(1-\tfrac{C}{\lambda})^{1/{d}}\Bigr]^{1/{d}}}.
\]
\end{prop}
\begin{proofb}
By Lemma \ref{Erlangfixedpoint}, we know that
\begin{equation}\label{eqretrial}
\pi_C =  Erl_C\bigl(\lambda\,g(\pi_C)\bigr),
    \quad g(x)=\sum\limits_{i=0}^{d-1}x^i,
\end{equation}
and by the mass--balance equation,
\begin{equation}\label{massr}
\lambda(1-\pi_C^d)=\sum_{k=1}^C\pi_k.
\end{equation}
From \eqref{massr}, since $ \pi_k \leq 1$ for all $1\leq k \leq C$, we have that $\pi_C^d\geq1-\tfrac{C}{\lambda}$ and therefore
\[
\pi_C\geq(1-\tfrac{C}{\lambda})^{1/d}.
\]
Also, since $ \pi_k \geq \pi_C$ for all $1\leq k \leq C$, we have that
\[
\lambda(1-\pi_C^d)\geq C \pi_C \geq C (1-\tfrac{C}{\lambda})^{1/d},
\]
and solving for $\pi_C$ we obtain $\pi_C\leq\Bigl[1-\frac{C}{\lambda}(1-\tfrac{C}{\lambda})^{1/d}\Bigr]^{1/d}=:u_0$. Then,
\[
\lambda g(\pi_C)\leq\lambda g(u_0)=\lambda(1-u_0^d)/(1-u_0)=\frac{C(1-\tfrac{C}{\lambda})^{1/d}}{1-\Bigl[1-\frac{C}{\lambda}(1-\tfrac{C}{\lambda})^{1/d}\Bigr]^{1/d}}=\mathcal{R}_d(\lambda,C).
\]
By monotonicity of $Erl_C(\cdot)$, \eqref{eqretrial} implies $\pi_C\le Erl_C\bigl(\mathcal{R}_{d}(\lambda,C)\bigr)$, and raising to the power $d$ yields \eqref{upperboundretrial}.
\end{proofb}

\begin{remark}
By an argument analogous to the previous proof, but in the opposite direction, we obtain a lower bound on $\pi_C^d$:
\[
\pi_C^d\geq \Big[ Erl_C\Big(\frac{C}{1-(1-\tfrac{C}{\lambda})^{1/d}} \Big)\Big]^d.
\]
\end{remark}
Analogous to Proposition \ref{propretrial}, we now derive an Erlang--type lower bound for the stationary blocking probability $q_C^{d}$ under power-of-$d$. Recall that we denote $\boldsymbol{q}$ (rather than $\boldsymbol{\pi}$) for the fixed point under this policy.

\begin{prop}\label{proppower}
Fix $d\geq 1$, $C\in \mathbb{N}$ and $\lambda> C$. Let $\boldsymbol{q}$ be the fixed point of equation \eqref{eqpower} associated with the power-of-$d$ scheme. Then
\begin{equation}\label{lowerboundpower}
q_C^d \;\ge\; [ Erl_C(\mathcal P_d(\lambda,C))]^d,
\end{equation}
with
\[
\mathcal P_d(\lambda,C) =
\frac{ C\,(1 - \tfrac{C}{\lambda})^{1/d}}{\bigl(\tfrac{C}{\lambda}(1 - \tfrac{C}{\lambda})^{1/d} + 1 - \tfrac{C}{\lambda}\bigr)^{1/d}- (1 - \tfrac{C}{\lambda})^{1/d}}.
\]
\end{prop}

To prove the previous proposition, we will use the following result.

\begin{lemma}\label{remark_power}
Fix $d\geq 1$, $C\in \mathbb{N}$ and $\lambda>0$. Let $\boldsymbol{q}$ be the fixed point of \eqref{eqpower} for power-of-$d$, and let $\lambda'>0$ be such that
\[
\lambda\sum\limits_{i=0}^{d-1}q_{k-1}^{d-1-i}q_k^i\leq\lambda'\text{ for all } k=1,\dots,C.
\]
Then $q_C\leq Erl_C(\lambda')$. Analogously,
\[
\lambda\sum\limits_{i=0}^{d-1}q_{k-1}^{d-1-i}q_k^i\geq \lambda'\text{ for all } k=1,\dots,C
\]
implies $q_C\geq Erl_C(\lambda')$.
\end{lemma}
\begin{proofb}
Since
\[
q_{k-1}^d-q_{k}^d=(q_{k-1}-q_{k})\sum\limits_{i=0}^{d-1}q_{k-1}^{d-1-i}q_k^i,
\]
applying the assumption in \eqref{eqpower}, it follows that
\begin{equation}\label{almostErlang}
k(q_{k}-q_{k+1})
  \le \lambda'\,(q_{k-1}-q_{k}),
  \quad k=1,\dots,C.
\end{equation}
Denoting $d_k=q_k-q_{k+1}$ for $k=0,\dots,C-1$ and $d_C=q_C$ (since $q_{C+1}=0$), equation \eqref{almostErlang} implies
\[
d_k \le \frac{\lambda'}{k}\,d_{k-1},
  \quad k=1,\dots,C.
\]
Let us prove $q_C \le  Erl_C(\lambda')$ by contradiction. If it were not true, then
\[
\frac{(\lambda')^C/C!}{\displaystyle\sum_{i=0}^{C}(\lambda')^i/i!}<q_C=d_C,
\]
and by backward induction on $k$ we would have
\[
\frac{(\lambda')^k/k!}{\displaystyle\sum_{i=0}^{C}(\lambda')^i/i!}<d_k,
  \quad k=0,1,\dots,C.
\]
Summing over $k$ from $0$ to $C$ gives $1 <\sum\limits_{k=0}^C d_k=q_0-q_{C+1}=1,$ which is a contradiction.
\end{proofb}

\begin{proofb}[Proof of Proposition \ref{proppower}]
By the previous lemma, it suffices to show that $\lambda\sum\limits_{i=0}^{d-1}q_{k-1}^{d-1-i}q_k^i$ is lower bounded by $\mathcal P_d(\lambda,C)$. Since
\begin{equation}\label{neweq12}
\lambda\sum\limits_{i=0}^{d-1}q_{k-1}^{d-1-i}q_k^i\geq\lambda\sum\limits_{i=0}^{d-1}(q_{C-1}^{d-1-i}q_C^{i})=\lambda q_{C-1}^{d-1}\sum\limits_{i=0}^{d-1} \big( \tfrac{q_C}{q_{C-1}}\big)^i,
\end{equation}
using $\sum_{i=0}^{d-1} r^i = (r^d-1)/(r-1)$ with $r=\tfrac{q_C}{q_{C-1}}<1$, we obtain that the expression in \eqref{neweq12} is lower bounded by
\begin{equation}\label{neweq13}
\lambda\frac{q_{C-1}^d-q_C^d}{q_{C-1}-q_C}.
\end{equation}
From the fixed point equation with $k=C$, we know that
\[
Cq_C=\lambda(q_{C-1}^d-q_{C}^d)\text{ and }q_{C-1}=\bigl(\tfrac{C}{\lambda} q_C + q_C^d\bigr)^{1/d},
\]
and $q_C> 0$. Therefore \eqref{neweq13} is equal to
\begin{equation}\label{neweq14}
\frac{C}{\bigl(\tfrac{C}{\lambda}q_C^{1-d} + 1\bigr)^{1/d}-1}.
\end{equation}
Since $q_C^d\geq 1-\tfrac{C}{\lambda}$ (which follows from the mass--balance equation), it follows that
\[
q_C^{1-d}\leq (1-\tfrac{C}{\lambda})^{\frac{1-d}{d}}
\]
and therefore, the expression in \eqref{neweq14} is greater than or equal to
\[
\frac{C}{\bigl(\frac{C}{\lambda}  (1-\tfrac{C}{\lambda})^{\frac{1-d}{d}}+ 1\bigr)^{1/d}-1},
\]
which can be rewritten as
\[
\frac{
 C\,(1 - \tfrac{C}{\lambda})^{1/d}
}{
  \bigl(\tfrac{C}{\lambda}(1 - \tfrac{C}{\lambda})^{1/d} + 1 - \tfrac{C}{\lambda}\bigr)^{1/d}
  - (1 - \tfrac{C}{\lambda})^{1/d}
},
\]
which is equal to $\mathcal P_d(\lambda,C)$.
\end{proofb}

In order to prove Theorem \ref{theorem2}, we introduce two auxiliary lemmas, whose proofs are included in the Appendix \ref{appendix_comparison}. For notational convenience, define $t=\tfrac{\lambda}{C}>1$. Then the previous expressions $\mathcal{P}_d(\lambda,C)$ and $\mathcal{R}_d(\lambda,C)$ can be rewritten as $\mathcal{P}_d(t,C)$ and $\mathcal{R}_d(t,C)$, with
\[
\mathcal{R}_{d}(t,C)=\frac{C(1-\tfrac{1}{t})^{1/{d}}}{1-\Bigl[1-\tfrac{1}{t}(1-\tfrac{1}{t})^{1/{d}}\Bigr]^{1/{d}}}
\]
and
\[
\mathcal{P}_d(t,C)=\frac{  C\,(1 - \tfrac{1}{t})^{1/d}
}{  \bigl(\tfrac{1}{t}(1 - \tfrac{1}{t})^{1/d} + 1 - \tfrac{1}{t}\bigr)^{1/d}
  - (1 - \tfrac{1}{t})^{1/d}}.
\]

\begin{lemma}\label{auxiliarlemmaasympt10}
Fix $C\in\mathbb{N}$ and $d\geq 1$. As $t\to\infty$, the following hold:
\[
\mathcal{R}_{d}(t,C)
= C\, d\,t + \frac{C(1- d)}2 + O(t^{-1}),
\ \
\mathcal{P}_{d}(t,C)
= C\,d\,t + \frac{C(1-d)}2 + O(t^{-1}).
\]
\end{lemma}

\begin{lemma}\label{auxiliarlemmaasympt}
Fix $C\in\mathbb{N}$. As $t\rightarrow +\infty$, if
\[
x = a_1\,t + a_2 + O(t^{-1})
\]
with $a_1>0$, then
\[
Erl_C(x)
= 1 - \frac{C}{a_1\,t} + \frac{C\,(a_2+1)}{a_1^2\,t^2} + O(t^{-3}).
\]
\end{lemma}

\begin{proofb}[Proof of Theorem \ref{theorem2}]
From Propositions \ref{propretrial} and \ref{proppower}, defining $t=\tfrac{\lambda}{C}>1$, it suffices to show that, given $\tilde d>d$, there exists $t_0=t_0(C,d,\tilde d)$ such that
\begin{equation}\label{des1}
[ Erl_C\bigl(\mathcal{R}_{\tilde{d}}(t,C)\bigr)]^{\tilde{d}}<[ Erl_C\bigl(\mathcal{P}_d(t,C)\bigr)]^d \ \forall t\geq t_0.
\end{equation}
By Lemma \ref{auxiliarlemmaasympt10},
\[
\mathcal{R}_{\tilde d}(t,C)
= C\, \tilde d\,t + \frac{C(1- \tilde d)}2 + O(t^{-1}),
\qquad \mathcal{P}_{d}(t,C) = C\,d\,t + \frac{C(1-d)}2 + O(t^{-1}).
\]
Applying Lemma \ref{auxiliarlemmaasympt} to the previous expressions, it follows
\begin{align*}
 Erl_C\bigl(\mathcal{R}_{\tilde d}(t,C)\bigr)
&= 1 - \frac{1}{\tilde d\,t} + \frac{C(1-\tilde d)+2}{2C\,\tilde d^2}\,t^{-2} + O(t^{-3}),\\
 Erl_C(\mathcal{P}_d(t,C))
&= 1 - \frac{1}{d\,t} + \frac{C(1-d)+2}{2C\,d^2}\,t^{-2} + O(t^{-3}).
\end{align*}
In addition, by the binomial theorem
\[
(1 + y)^n = \sum_{k=0}^n \binom{n}{k} y^k,
\]
for any constants $a_1, a_2$ we have
\[
\bigl[1 - \tfrac{a_1}{t} + \tfrac{a_2}{t^2} + O(t^{-3})\bigr]^n
= 1 - \frac{n\,a_1}{t}  + \frac{n\,a_2 + \tfrac{n(n-1)}{2}a_1^2}{t^2}
  + O(t^{-3}),
\]
and therefore,
\begin{align*}
\bigl[ Erl_C\bigl(\mathcal{R}_{\tilde d}(t,C)\bigr)\bigr]^{\tilde d}
&= 1 - \frac{1}{t} + \frac{1}{C\,\tilde d}\,t^{-2} + O(t^{-3}),\\
\bigl[ Erl_C(\mathcal{P}_{d}(t,C))\bigr]^{d}
&= 1 - \frac{1}{t} + \frac{1}{C\,d}\,t^{-2} + O(t^{-3}).
\end{align*}
Taking their difference we obtain
\[
\Delta(t)
= \bigl[ Erl_C(\mathcal{P}_{d}(t,C))\bigr]^{d}
  - \bigl[ Erl_C\bigl(\mathcal{R}_{\tilde d}(t,C)\bigr)\bigr]^{\tilde d}
  = \frac{\tilde d - d}{C\,d\,\tilde d}\,t^{-2}
  + O(t^{-3}),
\]
where the leading coefficient $\frac{\tilde d - d}{C\,d\,\tilde d}$ is positive (because $\tilde d>d$). Therefore, there exists $t_0$ such that
\[
\bigl[ Erl_C(\mathcal{R}_{\tilde d}(t,C))\bigr]^{\tilde d}
< \bigl[ Erl_C(\mathcal{P}_{ d}(t,C))\bigr]^{d}
\quad\text{for all }t>t_0.
\]
Indeed, if $\lvert O(t^{-3})\rvert\le \tfrac{K}{t^3}$ for $t\ge T$, then for $t\ge T$
\[
\Delta(t)
\ge \frac{\tilde d-d}{C\,d\,\tilde d}\,t^{-2}
  - \frac{K}{t^3}
= \frac{1}{t^3}\Bigl(\frac{\tilde d-d}{C\,d\,\tilde d}\,t - K\Bigr),
\]
so choosing
\[
t_0 = \max\!\Bigl\{T,\;\frac{K\,C\,d\,\tilde d}{\tilde d - d}\Bigr\}
\]
and $\lambda_0=t_0 C$ we obtain the claimed inequality.
\end{proofb}

\subsection{Comparison for \texorpdfstring{$C=2$ and $d=2$}{C=2 and d=2}}

In the previous subsection we prove that, for $\tilde d>d$, the stationary blocking probability in retrial-of-$\tilde d$ is strictly lower than the stationary blocking probability in power-of-$d$ for $\lambda$ sufficiently large. In this subsection we prove the following proposition, which asserts that in the case $C=2$, $d=2$ and $\tilde d =3$, (or equivalently, $\tilde d\geq 3$,) the same behavior holds for all $\lambda>0$.

\begin{prop}\label{prop_all_lambda}
Let $C=2$ and $\lambda>0$. Then, the stationary blocking probability under retrial-of-$3$ is strictly smaller than under power-of-$2$.
\end{prop}

\begin{proofb}
Fix $C=2$ and $\lambda>0$. In this case, the fixed-point equations \eqref{fixedpointretrial} for retrial-of-$3$ reduce to the following
\begin{equation}\label{retrial_fixed_C2}
\begin{aligned}
\pi_1-\pi_2&=\lambda(1-\pi_1)\,g(\pi_2),\\
2\pi_2&=\lambda(\pi_1-\pi_2)\,g(\pi_2),
\end{aligned}
\end{equation}
with $\pi_0=1$, $\pi_3=0$ and $g(u)=1+u+u^2$. Similarly, the fixed-point equations \eqref{eqpower} for power-of-$2$ reduce to
\begin{equation}\label{power_fixed_C2}
\begin{aligned}
q_1-q_2&=\lambda(1-q_1^2),\\
2q_2&=\lambda(q_1^2-q_2^2),
\end{aligned}
\end{equation}
with $q_0=1$ and $q_3=0$.

We first note that $\pi_2\in(0,1)$.  Indeed, if $\pi_2=0$, then the second equation in \eqref{retrial_fixed_C2} gives
\[
0=\lambda(\pi_1-\pi_2)g(\pi_2)=\lambda\pi_1.
\]
Since $\lambda>0$, we get $\pi_1=0$, and substituting into the first equation yields $0=\lambda$, a contradiction. If $\pi_2=1$, then $\pi_1=1$, and the second equation becomes
\[
2=\lambda(\pi_1-\pi_2)g(\pi_2)=0,
\]
again a contradiction. Hence $\pi_2\in(0,1)$.  Analogously, $q_2\in(0,1)$.

We now express $\lambda$ in terms of $\pi_2$ for retrial-of-$3$. From the first equation of \eqref{retrial_fixed_C2},
\[
\pi_1=\frac{\pi_2+\lambda g(\pi_2)}{1+\lambda g(\pi_2)}.
\]
Substituting into the second equation of \eqref{retrial_fixed_C2} gives
\begin{equation}\label{retrial_C2_4}
2\pi_2\bigl(1+\lambda g(\pi_2)\bigr)=\lambda^2 g(\pi_2)^2(1-\pi_2).
\end{equation}
Solving \eqref{retrial_C2_4} for $\lambda>0$ yields
\begin{equation}\label{lambda_from_pi2}
\lambda=\frac{\pi_2+\sqrt{\pi_2(2-\pi_2)}}{(1-\pi_2^3)}.
\end{equation}

On the other hand, from the second equation of \eqref{power_fixed_C2} we obtain the identity
\[
q_1^2=q_2^2+\frac{2q_2}{\lambda},
\]
which together with the first equation of \eqref{power_fixed_C2} gives
\begin{equation}\label{C2power_2}
\bigl(\lambda(1-q_2^2)-q_2\bigr)^2=q_2^2+\frac{2q_2}{\lambda}.
\end{equation}

We assume now that there exists $\lambda>0$ such that $\pi_2^3=q_2^2$, and we derive a contradiction. Since $\pi_2,q_2\in(0,1)$, there exists $s\in(0,1)$ such that $\pi_2=s^2$ and $q_2=s^3$. Plugging $\pi_2=s^2$ into \eqref{lambda_from_pi2} and using $\sqrt{s^2(2-s^2)}=s\sqrt{2-s^2}$, we obtain
\begin{equation}\label{lambda_s}
\lambda=\frac{s\bigl(s+\sqrt{2-s^2}\bigr)}{1-s^6}.
\end{equation}
Plugging $q_2=s^3$ and \eqref{lambda_s} into \eqref{C2power_2} and simplifying yields
\begin{equation}\label{equation_equal_s}
\bigl(s+r\bigr)^2\Bigl(s+r-2s^2\Bigr)=2(1-s^6),
\end{equation}
with $r=\sqrt{2-s^2}$, which is equivalent to
\begin{equation}\label{C2power_2.}
2(s^6-s^3-2s^2+3s-1)=-2(1+s^2-2s^3)r.
\end{equation}
Now let us show that \eqref{C2power_2.} cannot happen for any $s\in(0,1)$. To show this, it is sufficient to show that the polynomial
\[
D(s)=\bigl[2(s^6-s^3-2s^2+3s-1)\bigr]^2-\bigl[2(1+s^2-2s^3)\bigr]^2(2-s^2)
\]
obtained by subtracting the square of the right-hand side of \eqref{C2power_2.} from the square of the left-hand side of \eqref{C2power_2.}, has no root in $(0,1)$.

Expanding $D(s)$, we obtain
\begin{equation*}
D(s)=4(s^{12}-2s^{9}+2s^{7}-8s^{6}+8s^{5}-2s^{4}-2s^{3}+10s^{2}-6s-1),
\end{equation*}
which factorizes as
\[
D(s)=4(s-1)^2(s+1)Q(s),
\]
where
\[
Q(s)=s^9+s^8+2s^7+s^5+s^4-6s^3+2s^2-7s-1.
\]
Reordering the terms and extracting common factor, we can rewrite $Q(s)$ as
\[
Q(s)=s^3(s^6+s^5+2s^4+s^2+s-6)+2s(s-1)-5s-1,
\]
which is strictly negative for all $s\in(0,1)$. Then, since the other terms in the factorization of $D(s)$ are positive for all $s\in (0,1)$, we obtain $D(s)<0$ for all $s\in(0,1)$. Thus, \eqref{equation_equal_s} has no solution in $(0,1)$, contradicting $\pi_2^3=q_2^2$.

To finish the proof, let us consider the function
\[
T(\lambda)=q_2^2(\lambda)-\pi_2^3(\lambda),
\]
which gives, for each $\lambda$, the difference between the two stationary blocking probabilities (we write $q_2^2(\lambda)$ and $\pi_2^3(\lambda)$ to emphasize the dependence on $\lambda$). Since $T$ is continuous in $\lambda$ and $T(\lambda)\neq 0$ for all $\lambda>0$, it has constant sign on $(0,\infty)$. Therefore it is enough to show that for $\lambda$ sufficiently large, $T(\lambda)$ is positive, which follows from Theorem \ref{theorem2}.
\end{proofb}

\section{Numerics}
\label{section_numerics}
In this section we compare numerically the stationary blocking probabilities
\[
B_{\mathrm{retrial}}(\lambda)=\pi_C^{\tilde d}
\qquad\text{and}\qquad
B_{\mathrm{power}}(\lambda)=q_C^{d},
\]
corresponding to retrial-of-$\tilde d$ and power-of-$d$, respectively, in the regime $\tilde d>d$, with the same capacity $C$.

Recall that, by \eqref{fixedpointretrial} and \eqref{eqpower}, the stationary vectors
$\boldsymbol{\pi}=(\pi_k)_{k=0}^C$ and $\boldsymbol{q}=(q_k)_{k=0}^C$ are the unique solutions of
\[
k(\pi_k-\pi_{k+1})
=
\lambda(\pi_{k-1}-\pi_k)\Big(\sum_{j=0}^{\tilde d-1}\pi_C^{\,j}\Big),
\qquad 1\le k\le C,
\quad \pi_0=1,\ \pi_{C+1}=0,
\]
and
\[
k(q_k-q_{k+1})
=
\lambda(q_{k-1}^d-q_k^d),
\qquad 1\le k\le C,
\quad q_0=1,\ q_{C+1}=0,
\]
respectively.
Therefore, to compute $B_{\mathrm{retrial}}(\lambda)$ and $B_{\mathrm{power}}(\lambda)$ numerically,
it suffices to compute $\pi_C$ and $q_C$ from these analytic stationary equations.

To present the numerical procedure in a unified way, we denote by $\mathbf{P}=(P_k)_{k=0}^C$ the stationary tail vector,
where $P_k=\pi_k$ in the retrial model and $P_k=q_k$ in the power-of-$d$ model. In both cases, the stationary
equations can be rearranged to express $P_{k-1}$ in terms of $(P_k,P_{k+1})$. Thus, if $P_C$ were known, then together
with the boundary condition $P_{C+1}=0$ we could reconstruct $P_{C-1},P_{C-2},\dots,P_0$ by iterating backwards
from $k=C$ down to $k=1$.

\textbf{Retrial-of-$\tilde d$.} Given a candidate value $x\in(0,1)$ for $\pi_C$, we define $g_{\tilde d}(x)=\sum_{j=0}^{\tilde d-1} x^j.$
Setting $P_C=x$ and $P_{C+1}=0$, the recursion (obtained from \eqref{fixedpointretrial} with $d$ replaced by $\tilde d$) is
\begin{equation}\label{eq:backward-retrial}
P_{k-1}
=
P_k+\frac{k(P_k-P_{k+1})}{\lambda\,g_{\tilde d}(x)},
\qquad 1\le k\le C.
\end{equation}
Running \eqref{eq:backward-retrial} for $k=C,C-1,\dots,1$ yields an implied value $P_0(x)$.
The correct value must satisfy the boundary condition $P_0=1$, so $x=\pi_C$ is characterized as the solution of
\[
H_{\mathrm{retrial}}(x)=P_0(x)-1=0.
\]

\textbf{Power-of-$d$.} Given a candidate value $x\in(0,1)$ for $q_C$, set $P_C=x$ and $P_{C+1}=0$. From \eqref{eqpower} we obtain
\begin{equation}\label{eq:backward-power}
P_{k-1}^d
=
P_k^d+\frac{k(P_k-P_{k+1})}{\lambda},
\qquad 1\le k\le C,
\end{equation}
and hence $P_{k-1}$ is recovered by taking the nonnegative $d$-th root. Iterating \eqref{eq:backward-power}
backwards yields $P_0(x)$, and the correct value $x=q_C$ is characterized by
\[
H_{\mathrm{power}}(x)=P_0(x)-1=0.
\]

In both models, we solve $H(x)=0$ using Brent's method. 
The resulting values of $\pi_C$ and $q_C$ (and hence $B_{\mathrm{retrial}}(\lambda)$ and $B_{\mathrm{power}}(\lambda)$)
are interpreted as solutions of the analytic stationary equations up to solver tolerance.

We also compare the average probing cost of the two policies. Under power-of-$d$, each arriving job inspects $d$ servers in order to compare their occupancies, so the probing cost per arrival is exactly $d$. Under retrial-of-$\tilde d$, the number of probes is random, because the procedure stops as soon as a non-full server is found. Let $A_{\tilde d}$ denote the number of probes used by an arriving job under retrial-of-$\tilde d$, and let $p=\pi_C$ denote the stationary probability that a randomly selected server is full under this policy. The event $\{A_{\tilde d}\ge k\}$ means that the first $k-1$ probes have all selected full servers. Hence, for $1\le k\le \tilde d$,
\[
\mathbb P(A_{\tilde d}\ge k)=p^{k-1}.
\]
Using the tail-sum identity for positive integer-valued random variables,
\[
\mathbb E[A_{\tilde d}]
=
\sum_{k=1}^{\tilde d}\mathbb P(A_{\tilde d}\ge k)
=
\sum_{k=1}^{\tilde d}p^{k-1}
=
\sum_{k=0}^{\tilde d-1}p^k
=
\frac{1-p^{\tilde d}}{1-p},
\]
with the last expression understood as $\tilde d$ when $p=1$. Thus retrial-of-$\tilde d$ may have a worst-case probing cost larger than $d$, but its average probing cost can be smaller than $d$. This reflects the fact that, when a feasible server is found early, inspecting all sampled servers, as happens in power-of-$d$, may be unnecessary.

We define
\[
\phi_C(t)=B_{\mathrm{power}}(tC)-B_{\mathrm{retrial}}(tC),
\qquad t=\lambda/C.
\]
If a crossing exists, we define
\[
t_{\mathrm{crit}}(C)=\text{first root of }\phi_C(t)=0,
\qquad
\lambda_{\mathrm{crit}}(C)=C\,t_{\mathrm{crit}}(C).
\]
We emphasize $t_{\mathrm{crit}}(C)$ instead of $\lambda_{\mathrm{crit}}(C)$ because it is normalized by $C$ and remains directly comparable across capacities.

For each $C$ we evaluate $\phi_C(t)$ on a grid in $t$. If $\phi_C(t)$ changes sign on the grid,
we refine the root using Brent's method in $t$ to obtain $t_{\mathrm{crit}}(C)$.

\begin{figure}[t]
  \centering
  \includegraphics[width=0.98\linewidth]{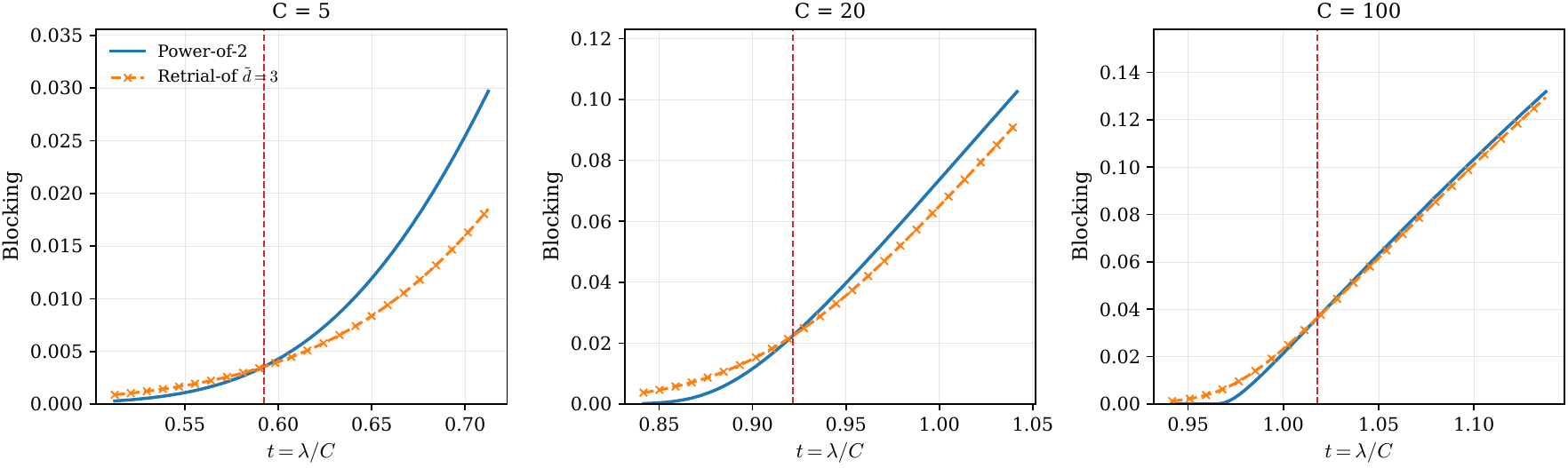}
  \caption{Blocking curves near the empirical crossing for $(d,\tilde d)=(2,3)$ and selected capacities.}
  \label{fig:blocking-vs-t}
\end{figure}

\begin{figure}[t]
  \centering
  \includegraphics[width=0.7\linewidth]{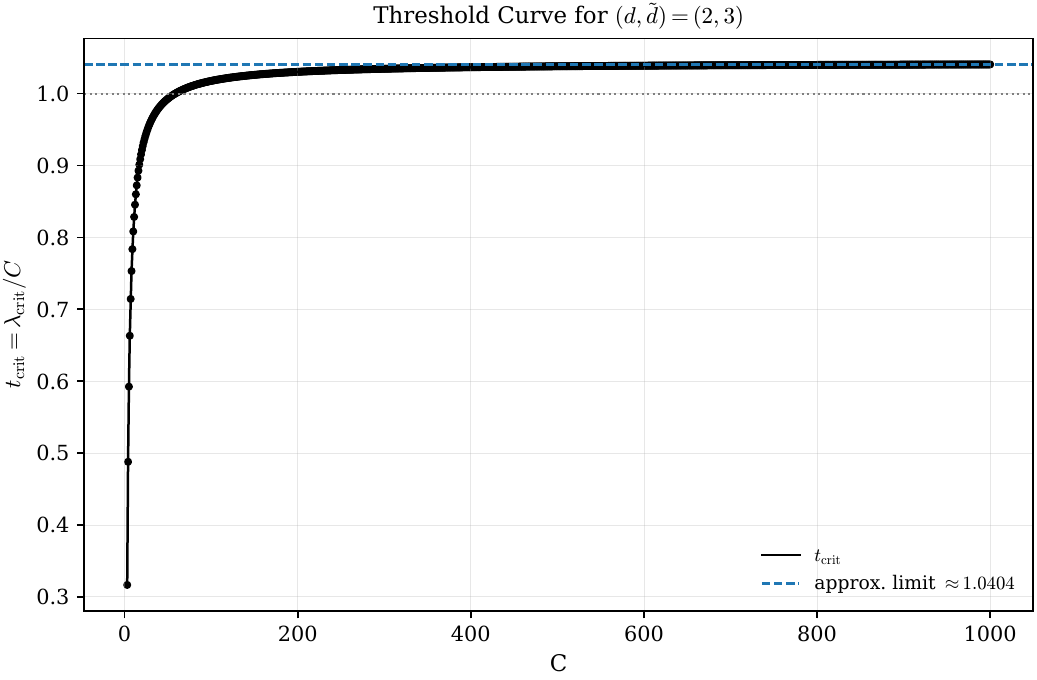}
  \caption{Empirical threshold $t_{\mathrm{crit}}(C)$ for $(d,\tilde d)=(2,3)$ and $C\geq 3$.}
  \label{fig:tcrit-vs-C}
\end{figure}

To quantify the improvement of retrial-of-$\tilde d$, relative to power-of-$d$, we report the percentage reduction
\[
\mathrm{reduction}_C(t)
=
100\,\frac{B_{\mathrm{power}}(tC)-B_{\mathrm{retrial}}(tC)}{B_{\mathrm{power}}(tC)}.
\]

For each row in Table~\ref{tab:reduction-max}, we define
\[
t^\ast=\arg\max_t \mathrm{reduction}_C(t),
\qquad
\lambda^\ast=C\,t^\ast,
\]
with maximization on the tested $t$-grid in the region where retrial is better (i.e.\ $t\ge t_{\mathrm{crit}}(C)$ when a crossing exists).

We also report the expected number of probes, $\mathbb E[A_{\tilde d}]$, used by retrial-of-$\tilde d$ at $\lambda^\ast$.
\begin{table}[t]
  \centering
  \begin{tabular}{ccccccccccc}
    \hline
    $C$ & $d$ & $\tilde d$ & $t_{\mathrm{crit}}$ & $\lambda_{\mathrm{crit}}$ & $t^\ast$ & $\lambda^\ast$ & $B_{\mathrm{power}}$ & $B_{\mathrm{retrial}}$ & reduction (\%) & $\mathbb E[A_{\tilde d}]$ \\
    \hline
    5 & 2 & 3 & 0.592 & 2.962 & 0.722 & 3.612 & 0.0334 & 0.0208 & 37.77 & 1.35 \\
    10 & 2 & 3 & 0.808 & 8.082 & 0.888 & 8.882 & 0.0490 & 0.0381 & 22.17 & 1.45\\
    20 & 2 & 3 & 0.922 & 18.432 & 0.982 & 19.632 & 0.0609 & 0.0535 & 12.11 & 1.52\\
    100 & 2 & 3 & 1.018 & 101.768 & 1.058 & 105.768 & 0.0696 & 0.0678 & 2.56 & 1.57\\
    5 & 2 & 4 & 0.444 & 2.220 & 0.624 & 3.120 & 0.0073 & 0.0011 & 84.88 & 1.22\\
    20 & 2 & 4 & 0.863 & 17.263 & 0.903 & 18.063 & 0.0130 & 0.0057 & 56.62 & 1.37\\
    100 & 2 & 4 & 0.977 & 97.694 & 0.987 & 98.694 & 0.0110 & 0.0085 & 22.49 & 1.42\\
    500 & 2 & 4 & 1.000 & 499.911 & 1.010 & 504.911 & 0.0144 & 0.0135 & 5.85 & 1.50\\
    10 & 2 & 5 & 0.672 & 6.719 & 0.782 & 7.819 & 0.0059 & 0.0004 & 93.25 & 1.26\\
    20 & 2 & 5 & 0.837 & 16.746 & 0.887 & 17.746 & 0.0067 & 0.0009 & 86.95 & 1.32\\
    100 & 2 & 5 & 0.969 & 96.857 & 0.979 & 97.857 & 0.0051 & 0.0019 & 62.60 & 1.40\\
    500 & 2 & 5 & 0.994 & 497.201 & 1.004 & 502.201 & 0.0093 & 0.0076 & 18.18 & 1.59\\
    \hline
  \end{tabular}
  \caption{Summary of thresholds, blocking values at $\lambda^\ast$, reduction maxima, and expected number of probes under retrial-of-$\tilde d$ at $\lambda^\ast$ for selected capacities and $(d,\tilde d)$ pairs.}
  \label{tab:reduction-max}
\end{table}

From Table~\ref{tab:reduction-max}, for fixed $(d,\tilde d)$ the improvement typically decreases with $C$ (for example, from $37.77\%$ at $C=5$ to $2.56\%$ at $C=100$ for $(d,\tilde d)=(2,3)$), while for fixed $d=2$ and fixed $C$, increasing $\tilde d$ lowers $t_{\mathrm{crit}}$ and increases the reduction (e.g., at $C=100$: $2.56\%$ for $\tilde d=3$, $22.49\%$ for $\tilde d=4$, and $62.60\%$ for $\tilde d=5$). In some parameter regimes the improvement is substantial; for example, for \(C=2\), \(d=2\), and \(\tilde d=5\), our numerical results indicate that the relative reduction in blocking can exceed \(90\%\).

The additional column $\mathbb E[A_{\tilde d}]$ shows that, in all reported cases, the expected number of probes under retrial-of-$\tilde d$ is smaller than $d=2$. Hence, at these operating points, retrial-of-$\tilde d$ not only achieves a smaller blocking probability than power-of-$2$, but also uses fewer probes on average, despite allowing a larger maximum number of probes.

In Figure~\ref{fig:blocking-vs-t}, the crossing between $B_{\mathrm{power}}(tC)$ and $B_{\mathrm{retrial}}(tC)$ is visible for representative capacities, and the gap after the crossing becomes smaller as $C$ increases. 

In Figure~\ref{fig:tcrit-vs-C}, the estimated threshold $t_{\mathrm{crit}}(C)$ for $(d,\tilde d)=(2,3)$ increases with $C$ and then seems to stabilize. More precisely, for $C=3,\dots,1000$ we found a root for every $C$, the sequence $t_{\mathrm{crit}}(C)$ is strictly increasing  and it appears to converge to approximately $1.0404$.

\section{Additional proofs}
\label{appendix}
\subsection{Proofs for Section \ref{section_uniqueness}}
\label{appendix_uniqueness}
This section of the Appendix collects the proof of Lemma \ref{Erlangfixedpoint} and the auxiliary Erlang--B identities used in the proof of Proposition \ref{lemmaA} (uniqueness of the fixed point).

\begin{proofb}[Proof of Lemma \ref{Erlangfixedpoint}]
Since $\boldsymbol{\pi}$ is a stationary point of the differential equation \eqref{eq.intro.4}, we have
\[
k(\pi_k-\pi_{k+1})
=\lambda(\pi_{k-1}-\pi_k)\sum_{j=0}^{d-1}\pi_C^{\,j},\qquad k=1,\dots,C,
\quad \pi_0=1,\ \pi_{C+1}=0.
\]
Let $d_k=\pi_k-\pi_{k+1}$ be the probability of having exactly $k$ jobs. Then
\[
k\,d_k = \Lambda\,d_{k-1}
\]
with $\Lambda= \lambda\sum\limits_{j=0}^{d-1}\pi_C^{\,j}$, so
\[
d_k=\frac{\Lambda^k}{k!}\,d_0.
\]
Since $\sum\limits_{k=0}^C d_k=1$, it follows that
\[
d_0=\Bigl(\sum\limits_{i=0}^{C}\frac{\Lambda^i}{i!}\Bigr)^{-1}
\]
and hence
\[
\pi_C=\tfrac{\Lambda^C/C!}{\sum_{i=0}^{C}\Lambda^i/i!}
= Erl_C(\Lambda).
\]
On the other hand, given $\pi_C$, the value of $\Lambda$ is fixed, so $(d_k)_{k=0}^C$ (and thus $\boldsymbol{\pi}$ via $\pi_k=\sum\limits_{i=k}^C d_i$) is univocally determined.
\end{proofb}

\begin{proofb}[Proof of Lemma \ref{lemmader}.]
Set $t_C(x)=\sum_{i=0}^C x^i/i!$. Since $Erl_C(x)=\dfrac{x^C/C!}{t_C(x)}$
and $t_C'(x)=t_{C}(x)-x^C/C!$ we have
\begin{align*}
Erl_C'(x)
&=\frac{C}{x}  Erl_C(x)-\frac{x^C/C!}{t_C(x)}  \frac{(t_{C}(x)-x^C/C!)}{t_C(x)}\\
&=\frac{C}{x} Erl_C(x)- Erl_C(x) (1- Erl_C(x))
= Erl_C(x)\Bigl[\frac{C}{x}-1+ Erl_C(x)\Bigr].
\end{align*}
Let us now prove that
\[
Erl_C(x)\bigl[C-x(1-Erl_C(x))\bigr]<1-Erl_C(x)
\]
by induction on $C$.

\textit{Case $C=1$:}
Since $Erl_1(x)=\frac{x}{1+x}$, we have
\[
Erl_1(x)\bigl[1-x(1-Erl_1(x))\bigr]
=\frac{x}{(1+x)^2}<\frac{1}{1+x}=1-Erl_1(x).
\]

\textit{Inductive hypothesis:}
Let us assume
\[
Erl_{C-1}(x)\bigl[C-1-x(1-Erl_{C-1}(x))\bigr]<1-Erl_{C-1}(x).
\]
We want to prove
\begin{equation}\label{indhyp4}
Erl_C(x)\bigl[C-x(1- Erl_C(x))\bigr]<1- Erl_C(x).
\end{equation}
Applying the recursive formula $Erl_C(x)=\frac{xErl_{C-1}(x)}{C+xErl_{C-1}(x)}$, it follows that \eqref{indhyp4} is equivalent to:
\begin{align}
\label{indhyp6}
\frac{xErl_{C-1}(x)}{C+xErl_{C-1}(x)}\bigl[C^2+CxErl_{C-1}(x)-xC\bigr]<C.
\end{align}
Multiplying both sides by $\frac{C+xErl_{C-1}(x)}{C}$, \eqref{indhyp6} can be rewritten as:
\begin{align}
\label{indhyp3}
xErl_{C-1}(x)\bigl[C-1-x(1-Erl_{C-1}(x))\bigr]<C.
\end{align}
By inductive hypothesis, the left term of \eqref{indhyp3} is lower than $x(1-Erl_{C-1}(x))$. So it is sufficient to prove that
\[
x(1-Erl_{C-1}(x))\leq C,
\]
which follows from the conservation law $x(1-Erl_{C-1}(x))\leq C-1$.
\end{proofb}

\subsection{Additional convergence proofs for Section \ref{section_convergence}}
\label{appendix_convergence}
We start this section by proving the Lipschitz property of the function $\boldsymbol{F}$ established in Lemma \ref{lemmaD}.

\begin{proofb}[Proof of Lemma \ref{lemmaD}.]
Let $x_{C+1}=y_{C+1}=0$ and $g(u)=\sum_{i=0}^{d-1}u^i$. For $\mathbf x,\mathbf y\in\mathcal S\subset[0,1]^{C+1}$,
\[
\begin{aligned}
|\mathbf F(\mathbf x)-\mathbf F(\mathbf y)|
&=\Big|\sum_{k=1}^{C}\lambda(x_{k-1}-x_k)g(x_C)-k(x_k-x_{k+1}) -\lambda(y_{k-1}-y_k)g(y_C)+k(y_k-y_{k+1})\Big| \\
&\le 2\sum_{k=1}^{C} k|x_k-y_k|+ \lambda\sum_{i=0}^{d-1}\sum_{k=1}^{C}\big|(x_{k-1}-x_k)x_C^{\,i}-(y_{k-1}-y_k)y_C^{\,i}\big|.
\end{aligned}
\]
Using $|ab-cd|\le|a-c||b|+|c||b-d|$, $|x_C^{\,i}|\le1$, and $|y_{k-1}-y_k|\le|y_{k-1}|+|y_k|\le2$, it follows that
\[
\begin{aligned}
|\mathbf F(\mathbf x)-\mathbf F(\mathbf y)|
\le\;& 2C|\mathbf x-\mathbf y|
+ \lambda\sum_{i=0}^{d-1}\sum_{k=1}^{C}\Big(|(x_{k-1}-y_{k-1})-(x_k-y_k)|+2|x_C^{\,i}-y_C^{\,i}|\Big)\\
\le\;& 2C|\mathbf x-\mathbf y|
+ 2\lambda d\,|\mathbf x-\mathbf y|
+ 2\lambda C \sum_{i=0}^{d-1}|x_C^{\,i}-y_C^{\,i}|.
\end{aligned}
\]
Since $|x_C^{\,i}-y_C^{\,i}|\le i\,|x_C-y_C|$ and $2\sum_{i=0}^{d-1}i=d(d-1)$, we obtain
\[
|\mathbf F(\mathbf x)-\mathbf F(\mathbf y)|
\le (2C+2\lambda d)|\mathbf x-\mathbf y|
+ \lambda d(d-1)C\,|x_C-y_C|
\le M|\mathbf x-\mathbf y|,
\]
with $M=2C+2\lambda d+\lambda d(d-1)C$.

The proof that $\mathbf F$ is locally Lipschitz in $\mathbb{R}^{C+1}$ is analogous:
given $\mathbf{x}\in\mathbb{R}^{C+1}$, consider the neighborhood
\[
K_R=\{\mathbf{y}\in\mathbb{R}^{C+1}:\ \|\mathbf{y}\|_\infty\le R\}
\]
with $R=1+\|\mathbf{x}\|_\infty$. Repeating the previous estimates with the uniform bound
$|y_i|\le R$ (in place of $|y_i|\le 1$), we obtain the desired result.
\end{proofb}

We now prove the results associated to convergence of solutions.

\begin{proofb}[Proof of Lemma \ref{lemmaB}.]
By Lemma \ref{lemmaD}, $\boldsymbol{F}$ is locally Lipschitz on an open neighborhood $V\subset\mathbb R^{C+1}$ of $\mathcal S$ (and this open set can be extended to $V\to\mathbb R^{C+1}$ considering $f:\mathbb R\times V\to\mathbb R^{C+1}$ with $f(t,x)=\boldsymbol{F}(x)$). Then, by the Picard--Lindel\"of theorem (Theorem 2.2 in \cite{Teschl}), for any initial condition $x(0)=x^0\in V$ there exists $\varepsilon>0$ and a unique local solution on a neighborhood $(-\varepsilon,\varepsilon)$ of $t=0$.
Then, by the maximal-solution theorem (Theorem 2.13 in \cite{Teschl}), there is a unique maximal solution $x(\cdot)$ defined on some interval $(T_-,T_+)$ with $0\in (T_-,T_+)$. Moreover, $\mathcal S$ is forward invariant: defining $g(x)=\sum_{j=0}^{d-1}x^j$, and denoting $s_k(t)$ and $\dot s_k(t)$ by $s_k$ and $\dot s_k$ respectively, from \eqref{eq.intro.4} it follows that:
\begin{itemize}
\item if $s_{k-1}=s_{k}$, then $\dot s_{k-1}-\dot s_k = \lambda (s_{k-2}-s_{k-1})  g(s_C) + k(s_k-s_{k+1})\ge0$;
\item if $s_k=0$, then $s_{k+1}=0$ and $\dot s_k=\lambda s_{k-1} g(s_C)\ge0$;
\item if $s_k=1$, then $s_{k-1}=1$ and $\dot s_k=-k(1-s_{k+1})\le0$.
\end{itemize}
Hence on every boundary face the vector field points inward, so the trajectory starting at $\mathbf s^0\in\mathcal S$ remains in $\mathcal S$. Since $\mathcal S$ is compact, and for each finite $T$ the set $[0,T]\times\mathcal S\subset \mathbb R\times V$ is compact, Corollary 2.16 in \cite{Teschl} implies that a maximal solution can not terminate in finite time; thus $T_+=+\infty$.

To prove uniqueness, suppose there exist two solutions $S(t)$ and $T(t)$ with $S(0) = T(0)$ and define
\[
t_0 = \inf \big\{ t > 0 : S(t) \neq T(t) \big\}.
\]
By continuity, $S(t_0) = T(t_0)$.
Applying Picard--Lindel\"of at $t_0$, there exists $\delta > 0$ such that $S(t) = T(t)$ on $(t_0- \delta, t_0 + \delta)$, a contradiction. Hence, the global solution is unique.
\end{proofb}

\medskip

\begin{proofb}[Proof of Lemma \ref{lemmaB2}]
For any initial values $\mathbf{s}^{0}\in\mathcal S$, define the two initial conditions $\mathbf{s}^{u}$ and $\mathbf{s}^{\ell}$ by
\[
s_{k}^{u}=\max\{s_{k}^{0},\pi_{k}\},\qquad
s_{k}^{\ell}=\min\{s_{k}^{0},\pi_{k}\},\qquad 0\le k\le C,
\]
and let $\mathbf{s}^{u}(t)$ and $\mathbf{s}^{\ell}(t)$ denote the solutions with initial conditions $\mathbf{s}^{u}$ and $\mathbf{s}^{\ell}$, respectively.
Since $s_{k}^{\ell}\le s_{k}^{0}\le s_{k}^{u}$ for $k=1,2,\ldots,C$, it follows by the same argument as in Lemma 3 of \cite{Srikant}, that the order is preserved for all $t$, i.e.,
\[
s_{k}^{\ell}(t)\le \pi_{k}\le s_{k}^{u}(t)\qquad\text{for all }t\text{ and all }k.
\]
Thus it is sufficient to show that
\[
\lim_{t\to\infty}\big|\mathbf{s}^{u}(t)-\boldsymbol{\pi}\big|
=\lim_{t\to\infty}\big|\mathbf{s}^{\ell}(t)-\boldsymbol{\pi}\big|=0,
\]
where $|\cdot|$ denotes the $\ell_{1}$ norm. This follows from the same argument as in Lemma 4 in \cite{Srikant}.
\end{proofb}

\medskip

\begin{proofb}[Proof of Lemma \ref{lemmaC}.]
For the case with replacement the argument is immediate from the law of large numbers for density--dependent Markov chains, Theorem 2.1 (chapter 11) in \cite{Ethier}, or Theorem 8.1 in \cite{Kurtz}. The process $\mathbf S^{N}(t)$ is a density dependent Markov chain with jumps of size $1/N$, all the rates are bounded by the constant $\lambda d+C$, and the drift $F$ is Lipschitz.

For the case without replacement, we repeat the same proof of Theorem 2.1 in \cite{Ethier} and only track the sampling correction. In the case with replacement, appears the factor $g(s_C)$ in the arrival rate, where $g(u)=\sum_{j=0}^{d-1}u^j.$ In the case without replacement, the only difference is that the function $g$ is replaced by $\tilde{g_N}$, with $\tilde{g_N}(u)=\sum\limits_{j=0}^{d-1}\frac{\binom{Nu}{j}}{\binom{N}{j}}\frac{N}{N-j}$. And for fixed $d$, we have the uniform bound
\[
\sup_{u\in[0,1]}\bigl|\tilde{g_N}(u)-g(u)\bigr|\;\le\;\frac{E_d}{N}
\]
for some positive constant $E_d$ which does not depend on $N$.
Consequently, the generator differs from the with replacement one by at most a constant in the rates, so the induced drift correction is of order $1/N$. In the time--change representation this adds only an $O(1)$ extension to the Poisson clocks (negligible under the Poisson strong law after division by $N$) and a deterministic error of size at most $A T/N$ on $[0,T]$. The martingale term still vanishes almost surely, and a Gr\"onwall estimate gives
\[
\sup_{t\le T}\bigl|{\mathbf S}^{N}(t)-\mathbf{s}\left(\mathbf{s}^{0}, t\right)\bigr|\;\longrightarrow\;0\quad\text{a.s.},
\]
with $\mathbf{s}\left(\mathbf{s}^{0}, t\right)$ the same fluid trajectory as in the case with replacement.
\end{proofb}

\begin{proofb}[Proof of Theorem \ref{theorem1}]
To show weak convergence of $\boldsymbol{\pi}^{(N)}$ to $\delta_{\boldsymbol{\pi}}$, we first note that the sequence $\{\boldsymbol{\pi}^{(N)}\}$ is tight due to the compactness of $\mathcal{S}$. By Prokhorov's theorem, every subsequence $\{\boldsymbol{\pi}^{(N_k)}\}$ contains a subsubsequence converging weakly to some limit $\bar{\boldsymbol{\pi}}$. By Skorokhod's representation theorem, we can couple random variables $S^{(N_k)}\sim{\boldsymbol{\pi}}^{(N_k)}$ and $S\sim{\boldsymbol{\bar\pi}}$ on a common probability space so that $S^{(N_k)}\to S$ almost surely.

Let $S^{(N_k)}(t)$ be the process started at $S^{(N_k)}(0)=S^{(N_k)}$, and let $S(t)$ be the solution of the limiting ODE with $S(0)=S$. By Lemma \ref{lemmaC}, convergence of the initial conditions implies convergence of the entire trajectories, so $S^{(N_k)}(t)\to S(t)$ weakly for each $t$. Since each ${\boldsymbol{\pi}}^{(N_k)}$ is stationary, we have $S^{(N_k)}(t)\sim {\boldsymbol{\pi}}^{(N_k)}$ for all $t$; letting $k\to\infty$, it follows that $S(t)\sim {\boldsymbol{\bar \pi}}$ for every $t$. Thus ${\boldsymbol{\bar \pi}}$ is a stationary distribution for the fluid-limit ODE, and by uniqueness of that stationary solution proved in Proposition \ref{lemmaA}, we conclude ${\boldsymbol{\bar \pi}}={\boldsymbol{\pi}}$.
\end{proofb}

\begin{proofb}[Proof of Proposition \ref{prop_chaos}]
Under ${\Pi}^{(N)}$, let $X_1^{(N)},\dots,X_N^{(N)}\in E$ denote the server occupancies, with $E=\{0,1,\dots,C\}$, a finite (hence Polish) space. By the symmetry of the routing policy considered, ${\Pi}^{(N)}$ is exchangeable, i.e., the joint law is invariant under any permutation of server indices. By Proposition 4.2 in \cite{Meleard:1996} for symmetric laws on $E^N$, propagation of chaos is equivalent to convergence of the empirical measures $\mu^{(N)}\in\mathcal P(E)$ to a deterministic limit, where
\[
\mu^{(N)}(k)=\frac{1}{N}\sum_{i=1}^N\mathbf 1_{\{X_i^{(N)}=k\}}.
\]
Since $\mu^{(N)}(k)=S_k^{(N)}-S_{k+1}^{(N)}$ for $0\le k\le C$
and, by Theorem \ref{theorem1}, the stationary law of $\mathbf S^{(N)}$ converges weakly to
$\delta_{\boldsymbol\pi}$, the continuous mapping theorem yields
\[
\mu^{(N)}\Rightarrow\rho
\quad\text{in }\mathcal P(E),
\]
where $\rho$ is the deterministic measure defined in \eqref{eqrho}. By this equivalence, ${\Pi}^{(N)}_{\mid l}\Rightarrow \rho^{\otimes l}$ on $E^p$.
\end{proofb}

\begin{proofb}[Proof of Proposition \ref{prop_limit_blocking}]
Recall that $\boldsymbol{\pi}^{(N)}$ is stationary distribution of the process $\mathbf{S}^{(N)}$. By Little's Law, the long-term average number of jobs in the system $L$ is equal to the effective arrival rate $\lambda_{\text{eff}}$ multiplied by the average time a job spends in the system $W$.
Here, $\lambda_{\text{eff}} = N\lambda(1-P_{block}^{(N)})$ (admitted jobs), and $W = 1/\mu = 1$ (mean service time, since we set $\mu=1$). Thus $L = N\lambda(1-P_{block}^{(N)})$. On the other hand, $L$ equals the expected number of jobs in the system, which is equivalent to $N\sum_{k=1}^C \pi_{k}^{(N)}$ (by the tail--sum identity: $\sum_{k=0}^C k\,\Pi^{(N)}_k=\sum_{i=1}^C \sum_{k=i}^C \Pi^{(N)}_k=\sum_{k=1}^C \pi^{(N)}_k.$) Equating both expressions we have
\begin{equation}
\label{eq.intro.3}
N \lambda\left(1-P_{block}^{(N)}\right)=N \sum_{k=1}^{C} \pi_{k}^{(N)}.
\end{equation}
Summing $F_{k}(\boldsymbol{\pi})=0$ over $k=1,\dots,C$ yields
\[
0=\sum_{k=1}^{C} F_{k}(\boldsymbol{\pi})
=\lambda\sum_{j=0}^{d-1}\pi_C^{j}(1-\pi_C)-\sum_{k=1}^{C} \pi_{k}
=\lambda(1-\pi_C^d)-\sum_{k=1}^{C} \pi_{k}.
\]
Thus $\lambda(1-\pi_C^{d})=\sum_{k=1}^{C} \pi_{k}$, whereas \eqref{eq.intro.3} gives $\lambda(1-P_{block}^{(N)})=\sum_{k=1}^{C}\pi_{k}^{(N)}$. By Theorem \ref{theorem1}, it follows that $\pi_{k}^{(N)}\to \pi_k$ for each $k$, therefore $P_{block}^{(N)}\to \pi_C^{\,d}$ as $N\to\infty$.
\end{proofb}

\subsection{Proof of the auxiliary lemmas related to Erlang used in Section \ref{section_comparison}}
\label{appendix_comparison}

\begin{proofb}[Proof of Lemma \ref{auxiliarlemmaasympt10}]
Let $x=\tfrac{1}{t}$. Let $A=1-x(1-x)^{1/d}$ and let
\[
M=C(1-x)^{1/ d}\quad\text{and}\quad
D=1-\bigl[1-x(1-x)^{1/d}\bigr]^{1/ d}=1-A^{1/ d}
\]
be the numerator and the denominator of $\mathcal{R}_{d}(t,C)$ respectively.
By the binomial expansion,
\[
(1-x)^{1/ d}
=1-\frac{x}{ d}+\frac{1- d}{2 d^2}\,x^2+O(x^3),
\]
so
\[
M=C\Bigl[1-\tfrac{x}{ d}+\tfrac{1- d}{2 d^2}\,x^2+O(x^3)\Bigr].
\]
Since $A=1-x\,(1-x)^{1/ d} =1-x+\tfrac{x^2}{ d}+O(x^3)$, we have
\[
A^{1/ d}=\bigl[1-(x-\tfrac{x^2}{ d}+O(x^3))\bigr]^{1/ d}
=1-\frac{x}{ d}+\frac{3- d}{2 d^2}\,x^2+O(x^3),
\]
and hence
\[
D=\frac{x}{ d}\Bigl[1-\tfrac{3- d}{2 d}\,x+O(x^2)\Bigr].
\]
Using $(1-U)^{-1}=1+U+O(U^{2})$ with $|U|<1$, gives
\[
\frac{1}{D}
=\frac{ d}{x}\Bigl[1+\tfrac{3- d}{2 d}\,x+O(x^2)\Bigr],
\]
and multiplying by $M$ yields
\[
\mathcal{R}_{ d}(t,C)
=C\,\frac{ d}{x}\Bigl[1
+\Bigl(\tfrac{3- d}{2 d}-\tfrac{1}{ d}\Bigr)\,x+O(x^2)\Bigr]
=C\, d\,t+\frac{C(1- d)}2+O(t^{-1}).
\]
For $\mathcal{P}_{d}(t,C)$, the numerator is the same that $\mathcal{R}_{d}(t,C)$. The denominator $E$ of $\mathcal{P}_{d}(t,C)$ can be written as $B^{1/d}-(1-x)^{1/d}$ with
\[
B=x\,(1-x)^{1/d}+(1-x)
=1-\frac{x^2}{d}+O(x^3).
\]
Since
\[
B^{1/d}=1-\frac{x^2}{d^2}+O(x^3),  \quad
(1-x)^{1/d}=1-\frac{x}{d}+\frac{1-d}{2d^2}\,x^2+O(x^3),
\]
we obtain
\[
E=\frac{x}{d}-\frac{3-d}{2d^2}\,x^2+O(x^3)
=\frac{x}{d}\Bigl[1-\frac{3-d}{2d}\,x+O(x^2)\Bigr],
\]
which has the same asymptotic form as the denominator of $\mathcal{R}_{ d}(t)$. Therefore
\[
\mathcal{P}_{d}(t,C)
=C\,d\,t+\frac{C(1-d)}2+O(t^{-1}).
\]
\end{proofb}

\begin{proofb}[Proof of Lemma \ref{auxiliarlemmaasympt}]
First, we observe that as $x\rightarrow +\infty$, we have
\[
Erl_C(x) \;=\; 1 - \frac{C}{x} + \frac{C}{x^2} + O(x^{-3}).
\]
Indeed, factorizing the denominator as
\[
\sum_{i=0}^{C} \frac{x^{i}}{i!}
= \frac{x^{C}}{C!}\,\Bigl(1 \;+\;\frac{C}{x} \;+\;\frac{C(C-1)}{x^{2}}
+ O\!\bigl(x^{-3}\bigr)\Bigr),
\]
we get
\begin{equation}\label{erlangasymp}
Erl_C(x)
=\frac{1}{1 + \frac{C}{x} + \frac{C(C-1)}{x^{2}} + O(x^{-3})}.
\end{equation}
Applying $(1+U)^{-1}=1-U+U^{2}+O(U^{3})$ with
$U=\tfrac{C}{x}+\tfrac{C(C-1)}{x^{2}}+O(x^{-3})<1$ we obtain
\[
Erl_C(x)=  1 - \frac{C}{x} + \frac{C}{x^{2}} + O(x^{-3}).
\]
Then, for $x = a_1\,t + a_2 + O(t^{-1})$, using again $(1+U)^{-1}=1-U+O(U^{2})$ with $|U|<1$ we obtain
\[
\frac{C}{x}
=\frac{C}{a_1\,t+a_2+O(t^{-1})}
=\frac{C}{a_1t}\,\frac{1}{1+\tfrac{a_2 t^{-1}}{a_1}+O(t^{-2})}
=\frac{C}{a_1t}-\frac{Ca_2}{a_1^2t^2}+O(t^{-3})
\]
and
\[
\frac{C}{x^2}
=\frac{1}{C}\Bigl(\frac{C}{a_1t}-\frac{Ca_2}{a_1^2t^2}+O(t^{-3})\Bigr)^2
=\frac{C}{a_1^2t^2}+O(t^{-3}).
\]
Replacing in \eqref{erlangasymp} yields the claimed formula.
\end{proofb}

\section*{Acknowledgement}
This work was supported by the ANR under the France 2030 program, grant “NF-NAI: ANR-22-PEFT-0003”.
  
\bibliography{bibliography}

\end{document}